\documentclass[12pt]{amsart}

\usepackage{amsfonts}
\usepackage{amssymb}
\usepackage{amstext}

\usepackage{graphicx}
\usepackage{tikz}
\usepackage{caption}
\usepackage[]{amscd}
\usepackage{hyperref}
\newtheorem{theorem}{Theorem}[section]
\newtheorem{lemma}{Lemma}[section]

\usepackage{multicol}
\theoremstyle{definition}

\theoremstyle{remark}
\newtheorem{remark}{Remark}[section]
\title[Isogenies between some $K3$ surfaces]{Isogenies between some $K3$ surfaces }

\usepackage{amsmath}

\usepackage{amsmath}
\begin{document}
\title[Isogenies between $K3$ surfaces]
{Isogenies between $K3$ surfaces of the Ap\'ery-Fermi pencil}

\author[M.-J. Bertin]{Marie Jos\'e Bertin}
\address{Marie Jos\'e Bertin \\ Institut de Math\'ematiques de Jussieu-Paris Rive Gauche \\ Sorbonne Universit\'e \\ 4 Place Jussieu, 75252 PARIS Cedex 05, France}
\email{marie-jose.bertin@imj-prg.fr}

\author[O. Lecacheux]{Odile Lecacheux}
\address{Odile Lecacheux  \\ Institut de Math\'ematiques de Jussieu-Paris Rive Gauche \\ Sorbonne Universit\'e \\ 4 Place Jussieu, 75252 PARIS Cedex 05, France}
\email{odile.lecacheux@imj-prg.fr}


\begin{abstract}

Elliptic fibrations of $K3$ surfaces belonging to the Ap\'ery-Fermi pencil ($Y_k$) may have $2$ or $3$-torsion sections defining on $(Y_k)$ automorphisms $\tau$ of order $2$ or $3$. First we consider $Y_{k}/\tau$ \ for some fibrations of the singular
$K3$ surface $Y_{10}$ in the case of two-torsion sections and obtain as for
the singular surface $Y_{2}$ either  the Kummer  surface associated to
$Y_{10}$ or $Y_{10}$ itself.  This last case is associated to the complex
multiplication on $Y_{10}$.  We prove also that for all the fibrations of
$Y_{2}$ with $3$-torsion sections $Y_{2}/\tau=Y_{10}.$ Results are different
for $Y_{10}$ where we can obtain for $Y_{10}/\tau$ one of the two surfaces with
transcendental lattice $[4 \quad 0\quad  18]$ or $\left[2 \quad 0 \quad 36\right]$. We explicit
also the link between $3$-isogeny on a fibration and base change on
other fibrations.

\end{abstract}


\maketitle

\section{Introduction}
In the Ap\'ery-Fermi pencil $Y_k$ with transcendental lattice $U\oplus<12>$ ($U$ being the hyperbolic lattice) and defined by the equation

\[(Y_k) \qquad X+\frac{1}{X}+Y+\frac{1}{Y}+Z+\frac{1}{Z}=k,\]
two singular $K3$ surfaces, namely $Y_2$ and $Y_{10}$, retain our attention due to the relation between their transcendental lattices 
 $T(Y_2)$ and $T(Y_{10})$ \cite{Be}, 
\[T(Y_2)=\begin{pmatrix}
	2 & 0\\
	0 & 4\\
\end{pmatrix} \qquad T(Y_{10})=\begin{pmatrix}
	6 &  0\\
	0 & 12\\
\end{pmatrix} =T(Y_2)[3].\]
Even more, in their Besan{\text{\c{c}}}on paper \cite{BL3}, section 5.2.1, Bertin and Lecacheux obtained Weierstrass equations of two extremal elliptic fibrations of $Y_{10}$ by $3$-isogenies from two extremal elliptic fibrations of $Y_2$.

The present paper intends to go deeper into differences and links between these two $K3$ surfaces, one of the difficulties being, as explained in \cite{BL3}, the fact that we know all the elliptic fibrations of $Y_2$ but only some of them concerning $Y_{10}$. It is divided in two parts; the first part concerns the $2$-isogenies of $Y_{10}$ and the second part is devoted to $3$-isogenies of $Y_2$ and $Y_{10}$. By $2$ or $3$- isogenies between two $K3$ surfaces, we always mean isogenies induced on an elliptic fibration by $2$-torsion or $3$-torsion sections which can be viewed as isogenies of degree $2$ or $3$ on elliptic curves $E_t$ over $\mathbb C(t)$.

In a previous paper \cite{BL2}, Bertin and Lecacheux classified all the involutions on $Y_k$ given by $2$-torsion sections, on one side {\it{Morrison-Nikulin involutions $\iota$}} such that the quotient surface $Y_k/\iota$ is the Kummer surface $K_k$ associated to $Y_k$ with transcendental lattice $T_{K_k}=T_{Y_k}[2]$ and on the other side all the other involutions $\tau$ such that the quotient surface $S_k:=Y_k/\tau$ is not a Kummer surface with transcendental lattice $T_{S_k}=<-2>\oplus <2> \oplus <6>$. They also proved that by specialisation to $Y_2$ {\it{Morrison-Nikulin involutions $\iota$}} remain {\it{Morrison-Nikulin involutions}} of $Y_2$ and that $S_2=Y_2$, the remaining not specialised involutions of $Y_2$ either being ``self-involutions'' meaning they preserve the elliptic fibration or give another fibration of $Y_2$ or being {\it{Morrison-Nikulin involutions}}.
In the present paper, considering ``self-involutions'' as complex multiplications over $E_t$, we shall prove in section $3$ the following theorem.

\begin{theorem}\label{thm:th1}
  \begin{enumerate}
  \item The {\it{Morrison-Nikulin involutions}} of $Y_k$ still remain by specialisation {\it{Morrison-Nikulin involutions}} of $Y_{10}$ and the specialised surface $S_{10}$ is the surface $Y_{10}$.
  \item As in the $Y_2$ case, among the non specialised involutions of $Y_{10}$ apart from ``self-involutions'' we find also a {\it{Morrison-Nikulin involution }}.
\item `` Self-involutions'' can be used to obtain the Mordell-Weil lattice of the fibration.
    \end{enumerate}

  \end{theorem}

  Also in the paper \cite{BL2}, Bertin and Lecacheux proved that the translation by a $2$-torsion section of the elliptic fibration $\#16$ of $Y_{10}$ inducing a symplectic involution  gives to $Y_{10}$ a Shioda-Inose structure i.e. defines a rational map $Y_{10}\rightarrow K_{10}$ where $K_{10}$ is a Kummer surface with transcendental lattice $T_{K_{10}}=T_{Y_{10}}[2]$. Hence the elliptic fibration $\#16$ induces an elliptic fibration $f$ of $K_{10}$.

  Thus, in section $4$, using the Kummer structure of $K_{10}$, product of two isogenous elliptic curves with complex multiplication we identify this fibration $f$ computing sections of its Mordell-Weil lattice. Finally the dual $2$-isogeny allows us to obtain the Mordell-Weil lattice of the fibration $\#16$ of $Y_{10}$.




 
  The second part is mostly devoted to $3$-isogenies. Given a $K3$ surface $S$ with transcendental lattice $T_S$, it is important to understand when the transcendental lattice of the $3$-isogenous surface is $T_S[3]$. This is understood from a generalisation of the following result concerning $2$-isogenies.
  Suppose given a $K3$ surface $S$ with a two-torsion elliptic fibration $f_t$ and Weierstrass equation
  \[(f_t) \qquad Y^2=X(X^2+a(t) X +b(t));\]
  suppose $X=h$ is a parameter for another fibration $(f_h)$ with two singular fibers $II^*$ at $h=0$ and $h=\infty$, then the two-isogeny on the fibration $(f_t)$ induces an elliptic fibration on another $K3$ surface named $K$. Using the formula for the dual two-isogeny one can prove that $K$ is the surface obtained by base change $h=u^2$ of the fibration $(f_h)$. It follows from Shioda's result \ref{th:sh} that the transcendental lattice of $K$ satisfies $T_K=T_S[2]$. This previous hypothesis is fulfilled if $(f_t)$ has a singular fiber of type $I_{12}^*$ (\cite{BL2} for fibration $\#4$).


  In other words, a sufficient condition on singular fibers of an elliptic fibration of a $K3$ surface $S$ with $2$-torsion section and transcendental lattice $T_S$ , to get by $2$-isogeny a $K3$ surface $\tilde{S}$ with transcendental lattice $T_{\tilde{S}}=T_S[2]$ is $S$ having an elliptic fibration with an $I_{12}^*$ fiber.
  This result can be generalised for $3$-isogenies in section $5$. In particular we prove the following theorem.

  \begin{theorem}\label{3isog}
Given a Weierstrass equation $(E_t)$ of a $K3$ surface $S$ with a $3$-torsion section given by the point $(0,0)$

\[(E_t) \qquad Y^2+A(t)XY+B(t)Y=X^3\]
and transcendental lattice $T_S$, denote $\tilde{S}$ the $3$-isogenous $K3$ surface $S/<\omega>$ where $\omega=(0,0)$ and $T_{\tilde{S}}$ its transcendental lattice.

A sufficient condition to get $T_{\tilde{S}}=T_S[3]$ is that $Y$ will be the elliptic parameter for a new elliptic fibration with two $II^*$ fibers at $0$ and $\infty$. In particular this condition can be realised if $(E_t)$ has singular fibers either $I_{18}$ or $I_{12}$ and $IV^*$.

\end{theorem}

With the help of this theorem or directly by computation when $T_{\tilde{S}}\neq T_S[3]$ we investigate the $3$-isogenies for fibrations of $Y_2$ and some fibrations of $Y_{10}$.

\begin{theorem}
  \begin{enumerate}
  \item The $3$-isogenous elliptic fibrations $H_{\#19}(k)$ (resp.$H_{\#20}(k)$) of the two generic fibrations $E_{\#19}(k)$ (resp. $E_{\#20}(k)$) are elliptic fibrations of the same $K3$ surface $N_k$ with transcendental lattice $U(3)\oplus <4>$, $U$ being the hyperbolic lattice.
  \item The specialised $K3$ surface $N_2$ is the $K3$ surface $Y_{10}$.
    \item The specialised $K3$ surface $N_{10}$ is the $K3$ surface with transcendental lattice $[4 \quad 0 \quad 18]$ (Shimada notation \cite{SZ}).
\end{enumerate}
  \end{theorem}

  \begin{theorem}
    The $K3$ surface $Y_k^{(3)}$ with transcendental lattice $U(3)\oplus <36>$ has a genus one fibration without section such that its Jacobian variety $J(Y_k^{(3)})=N_k$ is the $K3$ surface with transcendental lattice $U(3)\oplus <4>$.
    
\end{theorem}

\begin{theorem}
All the $3$-isogenous surfaces obtained from $3$-torsion sections of $Y_2$ are the $K3$-surface $Y_{10}$.

\end{theorem}

Applying Theorem 1.2 to fibration $21-c$ of $Y_2$ (see \cite{BL2} or section 5.1 of the present paper) satisfying $A(t)=t^2+5$ and $B(t)=1$, we obtain a fibration with two singular fibers $II^*$ giving by Shioda's theorem \ref{th:sh} and base change $h=u^3$ a rank $7$ elliptic fibration $E_u$ of $Y_{10}$ for which we can determine the Mordell-Weil lattice in section 5.2. This fibration $E_u$
\[ (E_u) \qquad Y^2=X^3-5u^2X^2+u^3(u^3+1)^2\]
is particularly interesting since, by a $2$-neighbour method (\cite{Kum},\cite{BE}), we get a two-torsion rank $4$ elliptic fibration $E_m$ of $Y_{10}$ (see section 5.4)
\[ (E_m) \qquad y^2=x^3-(m^3+5m^2-2)x^2+(m^3+1)^2x\]
giving the Morrison-Nikulin involution of Theorem 1.1 (2). 

   In particular, base changes in some elliptic fibrations of $Y_2$ lead to high rank elliptic fibrations of $Y_{10}$.

  \begin{theorem}

    \begin{enumerate}
      \item The elliptic fibration $E_h$ of $Y_2$ with Weierstrass equation 
      \[(E_h) \qquad Y^2=X^3-5h^2X^2+h^5(h+1)^2,\] 
with singular fibers $2II^*(0,\infty)$, $I_2(1)$, $2I_1$, gives by base change $h=u^3$ a rank $7$ elliptic fibration of $Y_{10}$ with singular fibers 
$2I_{0}^{\ast}\left(  \infty,0\right)$, $3I_{2}\left(  u^{3}+1\right)$
, $6I_{1}.$
    \item The elliptic fibration $E_f$ of $Y_2$ with Weierstrass equation
      \[(E_f) \qquad Y^2=X^3-f(2f-3)X^2+3f^2(f-1)^2X+f^3(f-1)^4,\]
      with singular fibers $II ^*(\infty)$, $III^*(0)$, $I_4(1)$, $I_1(32/27)$, gives by base change $f=u'^3$ a rank $4$ elliptic fibration of $Y_{10}$ with singular fibers $I_0^*(\infty)$, $III(0)$, $3I_4(1,\zeta_3,\zeta_3^2)$ ($\zeta_3^3=1$), $3I_1$.
    \item The elliptic fibration $E_g$ of $Y_2$ with Weierstrass equation
      \[ (E_g) \qquad y^2=x^3+4g^2x^2+g^3(g+1)^2x,\]
      with singular fibers $2III^*(0,\infty)$, $I_4(-1)$, $I_2(1)$, gives by base change $g=n^3$ a rank $4$ elliptic fibration of $Y_{10}$ with singular fibers $2III(0,\infty)$, $3I_4(-1,n^2-n+1)$, $3I_2(1,n^2+n+1)$.

      \end{enumerate}      

  \end{theorem}

The situation is totally different concerning $3$-isogenous surfaces of $Y_{10}$. We shall prove in section $6$ the following theorem.
\begin{theorem}
  \begin{enumerate}
  \item Elliptic fibrations numbered $80$ (See Shimada's notation \cite{SZ}) and $(11)$ \cite{BL3} induce by $3$-isogeny an elliptic fibration of a $K3$ surface with transcendental lattice $[4 \quad 0 \quad 18]$.

  \item There is an elliptic fibration of $Y_{10}$ with a $3$-torsion section whose $3$-isogenous surface is a $K3$ surface with transcendental lattice $[2 \quad 0 \quad 36]$.
\end{enumerate}

  \end{theorem}

Finally, in section $7$, we put our results on $2$ and $3$-isogenies on $Y_2$ and $Y_{10}$ in the perspective of a result of Boissi\`ere, Sarti and Veniani \cite{BSV} and find that all the $3$ isogenies between singular $K3$ surfaces involved in the previous sections are isometries of their correspondant rational transcendental lattices.

\bigskip

Most of the computations performed in this part aim at exhibiting new examples of $K3$ surfaces having a Shioda-Inose structure of order $3$ or being generalized Kummer surfaces of order $3$ \cite{GPM}.

The forecoming results will be published in another paper.

\section{Background}

We often use Shimada's notation \cite{SZ} for transcendental lattice i.e. $[a \quad b \quad c]$ meaning $\begin{pmatrix}
	a & b\\
	b & c\\
      \end{pmatrix}$.

The singular fibers of type $\,I_{n},D_{m},IV^{\ast}, ...$ at $t=t_{1},.,t_{m}$
or at  roots of a polynomial $p(t)$ of degree $m$ are denoted $mI_{n}%
(t_{1},..,t_{m})$ or $mI_{n}(p(t)).$ The zero
component of a reducible fiber is the component intersecting the zero section and is denoted 
$\theta_{0}$ or $\theta_{t_{0},0}.$ The other components denoted
$\theta_{t_{0},i}$ satisfy the property $\theta_{t_{0},i}\cdot\theta_{t_{0},i+1}=1.$

\subsection{Discriminant forms}

Let $L$ be a non-degenerate lattice. 
The \textbf{dual lattice} $L^*$ of $L$ is defined by
$$L^*:=\text{Hom}(L,\mathbb Z) =\{x\in L \otimes \mathbb Q /\,\,\, b(x,y)\in \mathbb Z \hbox{  for all }y \in L \}.$$
and the \textbf{discriminant group $G_L$} by
$$G_L:=L^*/L.$$
This group is finite if and only if $L$ is non-degenerate. In the latter case, its order is equal to the absolute value of the lattice determinant $\mid \det (G(e)) \mid$ for any basis $e$ of $L$.
A lattice $L$ is \textbf{unimodular} if $G_L$ is trivial.

Let $G_L$ be the discriminant group of a non-degenerate lattice $L$. The bilinear form on $L$ extends naturally to a $\mathbb Q$-valued symmetric bilinear form on $L^*$ and induces a symmetric bilinear form 
$$b_L: G_L \times G_L \rightarrow\mathbb Q / \mathbb Z.$$
If $L$ is even, then $b_L$ is the symmetric bilinear form associated to the quadratic form defined by
$$
\begin{matrix}
q_L: G_L  &\rightarrow  & \mathbb Q/2\mathbb Z\\
q_L(x+L) & \mapsto & x^2+2\mathbb Z.
\end{matrix}
$$
The latter means that $q_L(na)=n^2q_L(a)$ for all $n\in \mathbb Z$, $a\in G_L$ and $b_L(a,a')=\frac {1}{2}(q_L(a+a')-q_L(a)-q_L(a'))$, for all $a,a' \in G_L$, where $\frac {1}{2}:\mathbb Q/2 \mathbb Z \rightarrow \mathbb Q / \mathbb Z$ is the natural isomorphism.
The pair $\boldsymbol{(G_L,b_L)}$ (resp. $\boldsymbol{(G_L,q_L)}$) is called the \textbf{discriminant bilinear} (resp. \textbf{quadratic}) \textbf{form} of $L$.

When the even lattice $L$ is given by its Gram matrix, we can compute its discriminant form using the following lemma as explained in Shimada \cite{Shim1}.
\begin{lemma}\label{lem:gram}
Let $A$ the Gram matrix of $L$ and $U$, $V \in Gl_n(\mathbb Z)$ such that
\[ UAV=D=\begin{pmatrix}
        d_1 &  & 0\\
         & \ddots & \\
        0 &  & d_n\\
\end{pmatrix}
\]

with $1=d_1=\ldots =d_k < d_{k+1} \leq \ldots \leq d_n$. Then
\[G_L\simeq \oplus_{i>k} \mathbb Z/(d_i).\]
Moreover the $i$th row vector of $V^{-1}$, regarded as an element of $L^*$ with respect to the dual basis $e_1^*$, ..., $e_n^*$ generate the cyclic group $\mathbb Z/(d_i)$.
\end{lemma}

\subsection{Nikulin's results}


\begin{lemma}[Nikulin \cite{Nik}, Proposition 1.6.1]\label{L:nik0}
Let $L$ be an even unimodular lattice and $T$ a primitive sublattice. Then we have
$$G_T\simeq G_{T^{\perp}} \simeq L/(T\oplus T^{\perp}),\,\,\,\,\,\,\,q_{T^{\perp}}=-q_T.$$
In particular, $ |\det T| =|\det  T^{\perp}| =[L:T\oplus T^{\perp}]$.
\end{lemma}

\subsection{Shioda's results}\label{shio}
\begin{enumerate}
  \item We refer the reader to Shioda-Inose's paper \cite {SI}.

Let $\Phi : X \rightarrow \mathbb P^1$ be an elliptic surface with a section and consider the sections of this elliptic fibration  i.e. determined by the rational points of the corresponding elliptic curve defined over the field of rational functions in $s$.

Denote by $r(\Phi)$ the rank of the group of sections. Then the Picard number $\rho (X)$ satisfies the equation
\begin{equation*}
\rho(X) = r( \Phi)+2+\sum_{\nu =1}^h (m_{\nu}-1)
\end{equation*}
where $h$ is the number of singular fibers and $m_{\nu}$ the number of irreducible components of the corresponding singular fiber.

Besides, if $ r(\Phi)=0$ and if $n(\Phi)$ is the order of the group of sections of $\Phi$, one gets
\begin{equation*}
\mid \det T_X \mid = \mid \det S_X \mid =(\prod_
{\nu =1}^{h}m_{\nu}^{(1)})/n(\Phi )^2
\end{equation*}

where $m_{\nu}^{(1)}$ denotes the number of simple components of the divisor $D_{\nu}=\Phi ^{-1}(t_{\nu})$ of a singular fiber, $S_X$ the sublattice of algebraic cycles in $H_2(X,\mathbb Z)$ and $T_X=S_X^{\perp}$ the orthogonal complement of $S_X$ in $H_2(X,\mathbb Z)$  called the transcendental lattice.

\bigskip

\item Let $(S,\Phi, \mathbb P^1)$ be an elliptic surface with a section $\Phi$, without exceptional curves of first kind.

Denote by $NS(S)$ the group of algebraic equivalence classes of divisors of $S$\
.

Let $u$ be the generic point of $\mathbb P^1$ and $\Phi^{-1}(u)=E$ the elliptic\
 curve defined over $K=\mathbb C(u)$ with a $K$-rational point $o=o(u)$. Then, \
$E(K)$ is an abelian group of finite type provided that $j(E)$ is transcendental over $\mathbb C$.
Let $r$ be the rank of $E(K)$ and $s_1,..., s_r$ be generators of $E(K)$ modulo\
 torsion. Besides, the torsion group  $E(K)_{tors}$ is generated by at most two\
 elements $t_1$ of order $e_1$ and $t_2$ of order $e_2$ such that $1\leq e_2$, \
$e_2 |e_1$ and $\mid E(K) _{tors}\mid =e_1e_2$.

The group $E(K)$ of $K$-rational points of  $E$ is canonically identified with the group of sections of $S$ over $\mathbb P^1(\mathbb C)$.

For $s\in E(K)$, we denote by $(s)$ the curve image in $S$ of the section corresponding to $ s$.

Let us define
$$D_{\alpha}:=(s_{\alpha})-(o) \,\,\,\,\,\,1\leq\alpha \leq r$$
$$D'_{\beta}:=(t_{\beta})-(o) \,\,\,\,\,\,\beta =1,2.$$

Consider now the singular fibers of $S$ over
$\mathbb P^1$. We set
$$\Sigma:=\{v\in \mathbb P^1  / C_v=\Phi^{-1}(v)\,\,\,\, {\hbox {be a singular \                                                                                                                  
fiber}} \}$$
and for each $v\in \Sigma$, $\Theta_{v,i}$, $0\leq i \leq m_v-1$, the
$m_v$ irreducible components of $C_v$.

Let $\Theta_{v,0}$ be the unique component of $C_v$ passing through $o(v)$.

One gets
$$C_v=\Theta_{v,0} +\sum_{i\geq 1}\mu_{v,i}\Theta_{v,i},\,\,\,\,\,\,\,\mu_{v,i}\                                                                                                                  
\geq 1.$$
Let $A_v$ be the matrix of order $m_v-1$ whose entry of index
$(i,j)$ is $(\Theta_{v,i}\Theta_{v,j})$, $i,j\geq 1$, where $(DD')$
is the intersection number of the divisors $D$ et $D'$ along
$S$. Finally $f$ will denote a non singular fiber,
i.e. $f=C_{u_0}$ for $u_0\notin \Sigma$.
\begin{theorem}

The N\'eron-Severi group $NS(S)$ of the elliptic surface $S$
is generated by the following divisors
$$f, \Theta_{v,i} \,\,\,\,\,\,(1\leq i \leq                                                                                                                                                       
m_v-1,\,\,\,\,v\in\Sigma)$$
$$(o), D_{\alpha}\,\,\,\,\,\,1\leq \alpha \leq r,
\,\,\,\,D'_{\beta}\,\,\,\,\beta =1,2.$$                                                                                                                                                           
                                                                                                                                                                                                  
The only relations between these divisors are at most two relations                                                                                                                               
$$e_{\beta}D'_{\beta}\approx e_{\beta}(D'_{\beta} (o))f+\sum_{v\in
\Sigma
}(\Theta_{v,1},...,\Theta_{v,m_v-1})e_{\beta}A_v^{-1}\left ( \begin{array}{l}
(D'_{\beta}\Theta_{v,1})\\.\\
.\\
.\\
(D'_{\beta}\Theta_{v,m_v-1}
)
\end{array}
\right )
$$
where $\approx $ stands for the algebraic equivalence.
\end{theorem}
\end{enumerate}

\subsection{Shioda and Kuwata results}

\begin{theorem}\label{th:sh}\cite{Sh1}

Denote $F_{\alpha,\beta}^{(n)}$ the elliptic $K3$ surface defined by the Weierstrass equation
  \[F_{\alpha,\beta}^{(n)}\qquad y^2=x^3-3\alpha x+(t^n+\frac{1}{t^n}-2\beta) \qquad n=1,...,6\]
  $F_{\alpha,\beta}^{(2)} \simeq Km(E_1\times E_2)$ is the Shioda-Inose structure.

  The rank $r^{(n)}$ of $MW(F_{\alpha,\beta}^{(n)})$ satisfies

  \[r^{(n)}=h+Min(4(n-1),16)-\begin{cases} 0 & \text{if}\quad j_1\neq j_2\\
      n & \text{if} \quad j_1= j_2\neq 0,1\\
      2n & \text{if} \quad j_1= j_2=0 \ {\text{or}} \quad 1
    \end{cases}
  \]

  where $h=rk Hom(E_1,E_2)$.
  For all $\alpha,\beta$, for all $n\leq 6$,
  \[T(F_{\alpha,\beta}^{(n)})\simeq T(F_{\alpha,\beta}^{(1)})[n]\]
  \[\det(T(F_{\alpha,\beta}^{(n)}))=\det( T(F_{\alpha,\beta}^{(1)})).n^{\lambda}, \qquad \lambda=4-h\]

\end{theorem}

\begin{lemma}\cite{Ku1}
  The rank of the N\'eron Severi of $Km(E_1\times E_2)$ is $20$ if $E_1$ and $E_2$ are isogenous, and have complex multiplication. Moreover $h=2$ in that case.

  \end{lemma}

\subsection{Transcendental lattice}
Let $X$ be an algebraic $K3$ surface; the group $H^2(X,\mathbb Z)$, with the intersection pairing, has a structure of a lattice and by Poincar\'e duality is unimodular. The N\'eron-Severi lattice $NS(X):=H^2(X,\mathbb Z) \cap H^{1,1}(X)$ and the transcendental lattice $T(X)$, its orthogonal complement in $H^2(X,\mathbb Z)$ are primitive sublattices of $H^2(X,\mathbb Z)$ with respective signatures $(1,\rho -1)$ and $(2,20-\rho)$ where $\rho$ is the rank of the N\'eron-Severi lattice.

By Nikulin's lemma 2.2, their discriminant forms differ just by the sign, that is
\[(G_{T(X)},q_{T(X)})\equiv (G_{NS(X)}, -q_{NS(X)}).\]

\subsection{$2$-isogenous curves} \label{2-iso}
Let $E$ be an elliptic curve with a $2$-torsion point $\omega=(0,0)$
\[E:y^2=x^3+Ax^2+Bx, \quad B\neq 0\]
and $\phi$ the isogeny of kernel $\omega$. The $2$-isogenous elliptic curve $\phi(E)$ has a Weierstrass equation
\[Y^2=X^3-2AX^2+(A^2-4B)X\]
such that
\[(X=\frac{y^2}{x^2}, \quad Y=y\frac{B-x^2}{x^2}).\]

\subsection{$3$-isogenous curves}

\subsubsection{Method}

Let $E$ be an elliptic curve with a $3$-torsion point $\omega=\left(  0,0\right)                                                                                                                                      
$
\[E:Y^{2}+AYX+BY=X^{3}                             
\]
and $\phi$ the isogeny of \ kernel $<\omega>.$

To determine a Weierstrass \ equation for the elliptic curve $E/<\omega>$ we
need two functions $x_{1}$ of degree $2$ and $y_{1}$ of degree $3$ invariant
by $M\rightarrow M+\omega$ where $M=\left(  X_{M},Y_{M}\right)  $ \ is a
general point on $E.$ We compute $M+\omega$ and $M+2\omega\left(                                                                                                                                                 =M-\omega\right)  $ and can choose%

\begin{align*}
x_{1}  &  =X_{M}+X_{M+\omega}+X_{M+2\omega}=\frac{X^{3}+ABX+B^{2}}{X^{2}}\\                                                                                                                                       
y_{1}  &  =Y_{M}+Y_{M+\omega}+Y_{M+2\omega}
\\ &=\frac{Y\left(  X^{3}%
-AXB-2B^{2}\right)  -B\left(  X^{3}+A^{2}X^{2}+2AXB+B^{2}\right)  }{X^{3}}.%
\end{align*}
The relation between $x_{1}$ and $y_{1}$ gives a Weierstrass equation for
$E/<\omega>$
\[                                                                                                                                                                                                                 
y_{1}^{2}+\left(  Ax_{1}+3B\right)  y_{1}=x_{1}^{3}-6ABx_{1}-B\left(                                                                                                                                               
A^{3}+9B\right).                                     
\]
\subsubsection{Formulae} \label{formulae}

Denote $\zeta_3$ the $3$-th root of unity, $\zeta_3^{3}=1$.

Notice that the points with $x_{1}=-\frac{A^{2}}{3},y_{1}=\frac{%
-3}{2}B+\frac{1}{6}A^{3}\frac{\pm i\sqrt{3}}{18}\left( 27B-A^{3}\right)$ are $3$-torsion points.
Taking one of these points to origin and after some transformation we can
obtain a Weierstrass equation $y^{2}+ayx+by=x^{3}$ with the following transformations.

 Define
\begin{align*}                                                                                                                                            
S_{1}  &  =2\left(  \zeta_3^{2}-1\right)  y+6Ax-2\left(  \zeta_3-1\right)  \left(                                                                                                                                              
A^{3}-27B\right) \\                                                                                                                                                                                                
S_{2}  &  =2\left(  \zeta_3-1\right)  y+6Ax-2\left(  \zeta_3^{2}-1\right)  \left(                                                                                                                                              
A^{3}-27B\right)                                                                                                                                                                                                   
\end{align*}
and
\[                                                                                                                                                                                                                 
X=\frac{-1}{324}\frac{S_{1}S_{2}}{x^{2}},\qquad Y=\frac{1}{5832}\frac                                                                                                                                              
{S_{1}^{3}}{x^{3}}%
\]
then we get

\[                                                                                                                                                                                                                 
E/<\omega>:y^{2}+\left(  -3A\right)  yx+\left(  27B-A^{3}\right)  y=x^{3}.%
\]
If $A_{1}=-3A$, $B_{1}=27B-A^{3},$ then we define%

\begin{align*}                                                                                                                                                                                                     
\sigma_{1}  &  =2\left(  \zeta_3^{2}-1\right)  3^{6}Y+6A_{1}3^{4}X-2\left(                                                                                                                                               
              \zeta_3-1\right)  \left(  A_{1}^{3}-27B_{1}\right) \\                                                                                                                                                                    
\sigma_{2}  &  =2\left(  \zeta_3-1\right)  3^{6}Y+6A_{1}3^{4}X-2\left(                                                                                                                                                   
\zeta_3^{2}-1\right)  \left(  A_{1}^{3}-27B_{1}\right)                                                                                                                                                                   
\end{align*}
and then%

\[                                                                                                                                                                                                                 
x=\frac{-1}{324}\frac{\sigma_{1}\sigma_{2}}{3^{8}X^{2}}=-{\frac{3\,{X}^{3}%
+{A}^{2}{X}^{2}+3\,BAX+3\,{B}^{2}}{{X}^{2}}},\qquad y=\frac{1}{5832}%
\frac{\sigma_{1}^{3}}{3^{12}X^{3}}.%
\]

\subsubsection{Other properties of isogenies}\label{isog3}

The divisor of the function $Y$ is equal to $-3\left(  0\right)  +3\omega$ so
$Y=W^{3}$ where $W$ is a function on the curve $E/<\omega>$. If $X=WZ$, the
function field of $E/<\omega>$ is generated by $W$ and $Z.$ So replacing in
the equation of $E$ we obtain the relation between $Z$ and $W$

\[                                  
W^{3}+AZW+B-Z^{3}=0.                                                                                                                                             \]      
Let $\zeta_{3}=\frac{-1+\sqrt{-3}}{2}$ and $\zeta_{3}^{3}=1.$

From the equation
\[
W^{3}+AWZ+B=Z^{3}%
\]

with the transformation
\begin{align*}
W  & =-\frac{\sqrt{-3}}{9}\frac{y}{x}+\frac{A\zeta_{3}}{3}+\frac{\sqrt
{-3}\zeta_{3}^{2}}{9}\frac{\left(  27B-A^{3}\right)  }{x}\\
Z  & =-\frac{\sqrt{-3}}{9}\frac{y}{x}-\frac{A\zeta_{3}^{2}}{3}+\frac{\sqrt
{-3}\zeta_{3}}{9}\frac{\left(  27B-A^{3}\right)  }{x}%
\end{align*}
of inverse%
\begin{align*}
x  & =\frac{27B-A^{3}}{3W-3Z+A}\\
y  & =\frac{\zeta_{3}^{2}\left(  27B-A^{3}\right)  \left(  6Z-6\zeta_{3}%
^{2}W-2\zeta_{3}A\right)  }{2\left(  3W-3Z+A\right)  }%
\end{align*}
we recover the Weierstrass equation
\[
y^{2}-3Ayx+\left(  27B-A^{3}\right)  y=x^{3}.%
\]

\section{$2$-isogenies from elliptic fibrations of $Y_{10}$}\label{sec:1}
Let us prove the theorem 1.1.

\begin{theorem}\label{2iso}
  \begin{enumerate}
  \item The {\it{Morrison-Nikulin involutions}} of $Y_k$ still remain by specialisation {\it{Morrison-Nikulin involutions}} of $Y_{10}$ and the specialised surface $S_{10}$ is the surface $Y_{10}$.
  \item As in
    the $Y_2$ case, among the non specialised involutions of $Y_{10}$ apart from ``self-involutions'' we find also a {\it{Morrison-Nikulin involution }}.
  \item  `` Self-involutions'' can be used to obtain the Mordell-Weil lattice of the fibration.
    \end{enumerate}

 \end{theorem}

 \begin{proof}
   \begin{enumerate}
\item 
  With the same argument as for specialisation to $Y_2$, Morrison-Nikulin involutions specialised
to $Y_{10}$ remain Morrison-Nikulin involutions of $Y_{10}$.
Let us recall in the following Table 1 the Weierstrass equations of elliptic fibrations of $S_{10}$.

\begin{table}[tp]
\[%
\begin{array}
[c]{|c|c|}

\hline
\text{No} & \text{Weierstrass Equation} \\ \hline
\#7 &
\begin{array}
[c]{c}%
 E_7:{y}^{2}={x}^{3}+2{x}^{2}t\left (11t+1\right )-{t}^{2}\left (t-1                                              
\right )^{3}x \\ 

 III^*({\infty}), I_1^*(0), I_6(1), 2I_1(t^2+118t+25)\\

F_7:={Y}^{2}={X}^{3}-4{X}^{2}t\left (11t+1\right )+4{t}^{3}\left (118                                           
t+25+{t}^{2}\right )X \\
III^*                                                                               
                                                                                                     ({\infty}), I_2^*(0), I_3(1), 2I_2(t^2+118t+25)                                                                                                                                                                                                                                                                                                                                                                      
\end{array}
\\\hline\hline

 \#9 &

\begin{array}
[c]{c}%
 E_9:{y}^{2}={x}^{3}+28{t}^{2}{x}^{2}+{t}^{3}\left ({t}^{2}+98t+1                                                 
\right )x  \\         
2III^*(0,{\infty}), 2I_2(t^2+98t+1), 2I_1(t^2-98t+1)\\  
F_9:                                                                                                            
{Y}^{2}={X}^{3}-56{t}^{2}{X}^{2}-4{t}^{3}\left ({t}^{2}-98 t+1                                                  
\right )X    \\
  2III^*(0,{\infty}), 2I_2(t^2-98t+1), 2I_1(t^2+98t+1)\\
  \end{array}\\\hline\hline
\#14 &
\begin{array}
[c]{c}%

E_{14}:{y}^{2}={x}^{3}+t\left (98{t}^{2}+28t+1\right ){x}^{2}+{t}^{6}x

\\
 I_8^*(0), I_0^*(\infty), I_1(4t+1), I_1(24t+1), 2I_1(100t^2+28t+1)\\                                             
F_{14}:{Y}^{2}=X\left (X-96{t}^{3}-28{t}^{2}-t\right )\left (X-100{t}^{                                           
3}-28{t}^{2}-t\right)                                                                                           
\\

  I_4^*(0), I_0^*(\infty), I_2(4t+1), I_2(24t+1), 2I_2(100t^2+28t+1)\\

\end{array}
 \\\hline\hline
\#15 &
\begin{array}
[c]{c}%
E_{15}:{y}^{2}={x}^{3}-t\left (2+{t}^{2}-22t\right ){x}^{2}+{t}^{2}\left (t                                        
+1\right )^{2}x                                                                                                 
\\

I_1^*(0), I_4^*(\infty), I_4(-1), I_1(24), 2I_1(t^2-20t+4)\\
                                                                                                             
F_{15}:{Y}^{2}=X\left (X+t^3-24t^2\right ) \left (X+t^3-20t^2+4t    \right )                                                                                          
\\

I_2^*(0), I_2^*(\infty), I_2(-1), I_2(24), 2I_2(t^2-20t+4)\\

\end{array}\\\hline\hline
\#20 &
\begin{array}
[c]
{c}
E_{20}:{y}^{2}={x}^{3}+\left (\frac{1}{4}{t}^{4}-5{t}^{3}+{\frac {53}{2}}                                          
{t}^{2}-15t-\frac{3}{4}\right ){x}^{2}-t\left (t-10\right )x                                                    
\\

I_2(0), I_{12}(\infty), 2I_3(t^2-10t+1), I_2(10), I_1(1), I_1(9)\\

F_{20}:{Y}^{2}={X}^{3}+\left (-\frac{1}{2}{t}^{4}+10{t}^{3}-53{t}^{2}+30t+\frac{3}{2}                             
\right ){X}^{2}\\+\frac{1}{16}\left (t-1\right )\left (t-9\right )\left ({t}^{                                    
2}-10t+1\right )^{3}X \\

I_1(0), I_6(\infty), 2I_6(t^2-10t+1), I_2(1), I_2(9), I_1(10)\\

\end{array}

\\
\hline    
\end{array}
                                                                                                                                                                                             \]

 \caption{Fibrations $E_{i}$ of $Y_{10}$ and $F_{i}$ of  $S_{10}$: $S_{10}=Y_{10}$}\label{T:4}
 \end{table}
To prove that $S_{10}=Y_{10}$, we observe that the Weierstrass equation $E_9$ of an elliptic fibration of $Y_{10}$ has the same singular fibers as its 2-isogenous $F_9$. Moreover these two fibrations are isomorphic, the isomorphism being defined by $t=-T$, $x=-\frac{X}{2}$, $y=\frac{Y}{2\sqrt{-2}}$.
This property is sufficient to identify $S_{10}$ with $Y_{10}$.
\item
  There is a lot of `` self-isogenies'' preserving the elliptic fibration of $Y_{10}$. Let us give the rank $0$ ones.

 There are four $2$-isogenies from $Y_{10}$ to $Y_{10}$ defined by extremal elliptic fibrations with $2$-torsion sections, denoted number $8$, $87$, $153$, $262$ in Shimada-Zhang's paper \cite{SZ}. They are all " self-isogenies" preserving the fibration.

  We write below in Table 2, the Weierstrass equation $E$, its $2$-isogenous $F$ and the corresponding isomorphism.

  \begin{table}
    \centering
    \begin{tabular}{|c|c|}
   \hline

$E_{262}$ & $y^{2}=x^{3}+x^{2}(9(t+5)(t+3)+(t+9)^{2})-t^{3}(t+5)^{2}x$ \\ 
& $III^{\ast }(\infty )$, $I_{6}(0)$, $I_{4}(-5)$, $I_{3}(-9)$, $I_{2}(-4)$
\\ 
$F_{262}$ & $Y^{2}=X^{3}-2(9(T+5)(T+3)+(T+9)^{2})X^{2}+4(T+4)^{2}(T+9)^{3}X$
\\ 
& $III^{\ast }(\infty )$, $I_{6}(-9)$, $I_{4}(-4)$, $I_{3}(0)$, $I_{2}(-5)$
\\ 
& Isomorphism: $t=-T-9$, $x=-\frac{X}{2}$, $y=\frac{Y}{2\sqrt{-2}}$ \\ 
\hline
$E_{153}$ & $y^{2}=x^{3}+t(t^{2}+10t-2)x^{2}+(2t+1)^{3}t^{2}x$ \\ 
& $I_{3}(4)$, $I_{6}(-1/2)$, $I_{1}^{\ast }(0)$, $I_{2}^{\ast }(\infty )$ \\ 
$F_{153}$ & $Y^{2}=X^{3}-2T(T^{2}+10T-2)X^{2}+T^{3}(T-4)^{3}X$ \\ 
& $I_{6}(4)$, $I_{3}(-1/2)$, $I_{2}^{\ast }(0)$, $I_{1}^{\ast }(\infty )$ \\ 
& Isomorphism: $t=-\frac{2}{T}$, $x=-\frac{2X}{T^{4}}$, $y=-\frac{2\sqrt{-2}Y%
}{T^{6}}$ \\ 
\hline
$E_{87}$ & $y^{2}=x^{3}-(9t^{4}+9t^{3}+6t^{2}-6t+4)x^{2}+(21t^{2}-12t+4)x$
\\ 
& $I_{12}(\infty )$, $I_{6}(0)$, $2I_{2}(21t^{2}-12t+4)$, $%
2I_{1}(3t^{2}+6t+7)$ \\ 
$F_{87}$ & $Y^{2}=X^{3}-(9T^{4}+9T^{3}+6T^{2}-6T+4)X^{2}+(21T^{2}-12T+4)X$
\\ 
& $I_{12}(0)$, $I_{6}(\infty )$, $2I_{2}(3t^{2}+6t+7)$, $%
2I_{1}(21t^{2}-12t+4)$ \\ 
& Isomorphism: $t=-\frac{2}{T}$, $x=-\frac{2X}{9T^{4}}$, $y=\frac{2\sqrt{-2}%
}{27T^{6}}$ \\ 
\hline
$E_{8}$ & $y^{2}=x^{3}-(3t^{4}-60t^{2}-24)x^{2}-144(t^{2}-1)^{3}x$ \\ 
& $I_{2}(0)$, $2I_{3}(t^{2}+8)$, $I_{4}(\infty )$, $2I_{6}(t^{2}-1)$ \\ 
$F_{8}$ & $Y^{2}=X^{3}+2(3t^{4}-60t^{2}-24)X^{2}+9t^{2}(t^{2}+8)^{3}X$ \\ 
& $I_{2}(\infty )$, $2I_{3}(t^{2}-1)$, $I_{4}(0)$, $2I_{6}(t^{2}+8)$ \\ 
& Isomorphism: $t=\frac{2\sqrt{-2}}{T}$, $x=\frac{4X}{T^{4}}$, $y=\frac{8Y}{%
T^{6}}$.\\
  \hline
\end{tabular}
      \caption
      {Self Isogenies of extremal fibrations of $Y_{10}$}

      \end{table}


The rank $4$ elliptic fibration of $Y_{10}$  obtained in Theorem \ref{th:5.5} with singular fibers $3I_4$, $3I_2$, $2III$, Weierstrass equation 
\[(E_n) \qquad y^2=x^3+4n^2x^2+n(n^3+1)^2x\]
and its $2$-isogenous $E_n/\langle (0,0) \rangle$ are "self-isogenous", the elliptic fibration being preserved.

Indeed
\[E_n/\langle (0,0) \rangle \qquad Y^2=X^3-8T^2X^2-4T(T^3-1)^2X\]
with the same type of singular fibers. The isomorphism is given by
\[T=-n,\qquad Y=-2\sqrt{-2}y,\qquad X=-2x.\]

The following rank $4$ elliptic fibration  is an example of a non specialized Morrison-Nikulin involution of $Y_{10}$.
                                                                                                  The rank $4$ elliptic fibration  of $Y_{10}$ (\ref{rk4}), with Weierstrass equation
 \[(E_m) \qquad y^2=x^3-(m^3+5m^2-2)x^2+(m^3+1)^2x\]
 and singular fibers $I_0^*(\infty)$, $3I_4(m^3+1)$, $I_2(0)$, $4I_1(1,-5/3,m^2-4m-4)$, has a $2$-torsion section defining a Morrison-Nikulin involution from $Y_{10}$ to $K_{10}$, that is $F_m=E_m/\langle (0,0) \rangle$ is a rank $4$ elliptic fibration of $K_{10}$ with Weierstrass equation
 \[( F_m) \qquad Y^2=X^3+2(m^3+5m^2-2)X^2-(m^2-4m-4)(m-1)(3m+5)t^2X\]
 and singular fibers $I_0^*(\infty)$, $I_4(0)$, $7I_2(\pm 1,-5/3,m^2-m+1,m^2-4m-4)$.
 
Indeed, starting from $F_m$ and taking the new parameter $p=\frac{X}{(m^2-4m-4)(m-1)(3m+5)}$, we get a rank $1$ elliptic fibration with Weierstrass equation
\begin{align*}
(F_{p}) &  :Y^{2}=X^{3}+\frac{3}{4}p\left(  5p-1\right)  ^{2}X^{2}+\frac{1}%
{6}p^{2}\left(  2p-1\right)  \left(  5p-1\right)  \left(  49p^{2}%
-13p+1\right)  X\\
&  +\frac{1}{108}p^{3}\left(  2p-1\right)  ^{2}\left(  49p^{2}-13p+1\right)
^{2}
\end{align*}
and singular fibers $I_0^*(\infty)$, $I_3^*(0)$, $3I_3(\frac{1}{2}, 49p^2-13p+1)$.
The infinite section $P=(-\frac{1}{12} p(49p^2-13p+1), \frac{1}{8} p^2(49p^2-13p+1))$ is of height $h(P)=\frac{2}{3}$, is not equal to $2Q$ or $3Q$, hence the discriminant of the N\'eron-Severi lattice satisfies $\Delta=4 \times 4 \times 3^3 \times \frac{2}{3}=72 \times 4$.

Now we are going to compute its Gram matrix and deduce its transcendental lattice.

To compute the N\'eron-Severi lattice we order the generators as $(0)$ the zero section, $(f)$ the generic fiber, $\theta_i, 1 \leq i \leq 4$, $\eta_i, 1\leq i \leq 7$, $\gamma_i$, $\delta_i$, $\epsilon_i$, $1\leq i \leq 2$ the rational components of respectively $I_0^*$, $I_3^*$, the three $I_3$ and finally the infinite section $(P)$.
\[NS=\left ( \begin{smallmatrix}
	-2 & 1 & 0 & 0 & 0 & 0 & 0 & 0 & 0 & 0 & 0 & 0 & 0 & 0 & 0 & 0 & 0 & 0 & 0 & 0\\
	1 & 0 & 0 & 0 & 0 & 0 & 0 & 0 & 0 & 0 & 0 & 0 & 0 & 0 & 0 & 0 & 0 & 0 & 0 & 1\\
	0 & 0 & -2 & 0 & 1 & 0 & 0 & 0 & 0 & 0 & 0 & 0 & 0 & 0 & 0 & 0 & 0 & 0 & 0 & 0\\
	0 & 0 & 0 & -2 & 1 & 0 & 0 & 0 & 0 & 0 & 0 & 0 & 0 & 0 & 0 & 0 & 0 & 0 & 0 & 0\\
	0 & 0 & 1 & 1 & -2 & 1 & 0 & 0 & 0 & 0 & 0 & 0 & 0 & 0 & 0 & 0 & 0 & 0 & 0 & 0\\
	0 & 0 & 0 & 0 & 1 & -2 & 0 & 0 & 0 & 0 & 0 & 0 & 0 & 0 & 0 & 0 & 0 & 0 & 0 & 1\\
	0 & 0 & 0 & 0 & 0 & 0 & -2 & 0 & 1 & 0 & 0 & 0 & 0 & 0 & 0 & 0 & 0 & 0 & 0 & 0\\
	0 & 0 & 0 & 0 & 0 & 0 & 0 & -2 & 1 & 0 & 0 & 0 & 0 & 0 & 0 & 0 & 0 & 0 & 0 & 0\\
	0 & 0 & 0 & 0 & 0 & 0 & 1 & 1 & -2 & 1 & 0 & 0 & 0 & 0 & 0 & 0 & 0 & 0 & 0 & 0\\
	0 & 0 & 0 & 0 & 0 & 0 & 0 & 0 & 1 & -2 & 1 & 0 & 0 & 0 & 0 & 0 & 0 & 0 & 0 & 0\\
	0 & 0 & 0 & 0 & 0 & 0 & 0 & 0 & 0 & 1 & -2 & 1 & 0 & 0 & 0 & 0 & 0 & 0 & 0 & 0\\
	0 & 0 & 0 & 0 & 0 & 0 & 0 & 0 & 0 & 0 & 1 & -2 & 1 & 0 & 0 & 0 & 0 & 0 & 0 & 0\\
	0 & 0 & 0 & 0 & 0 & 0 & 0 & 0 & 0 & 0 & 0 & 1 & -2 & 0 & 0 & 0 & 0 & 0 & 0 & 1\\
	0 & 0 & 0 & 0 & 0 & 0 & 0 & 0 & 0 & 0 & 0 & 0 & 0 & -2 & 1 & 0 & 0 & 0 & 0 & 0\\
	0 & 0 & 0 & 0 & 0 & 0 & 0 & 0 & 0 & 0 & 0 & 0 & 0 & 1 & -2 & 0 & 0 & 0 & 0 & 0\\
	0 & 0 & 0 & 0 & 0 & 0 & 0 & 0 & 0 & 0 & 0 & 0 & 0 & 0 & 0 & -2 & 1 & 0 & 0 & 1\\
	0 & 0 & 0 & 0 & 0 & 0 & 0 & 0 & 0 & 0 & 0 & 0 & 0 & 0 & 0 & 1 & -2 & 0 & 0 & 0\\
	0 & 0 & 0 & 0 & 0 & 0 & 0 & 0 & 0 & 0 & 0 & 0 & 0 & 0 & 0 & 0 & 0 & -2 & 1 & 1\\
	0 & 0 & 0 & 0 & 0 & 0 & 0 & 0 & 0 & 0 & 0 & 0 & 0 & 0 & 0 & 0 & 0 & 1 & -2 & 0\\
	0 & 1 & 0 & 0 & 0 & 1 & 0 & 0 & 0 & 0 & 0 & 0 & 1 & 0 & 0 & 1 & 0 & 1 & 0 & -2\\
\end{smallmatrix} \right ).
 \]

Fom lemma \ref{lem:gram} we deduce that $G_{NS}=\mathbb Z/12 \mathbb Z \oplus \mathbb Z/ 24 \mathbb Z$ and we get generators $f_1$ and $f_2$ with respective norms $q(f_1)=-\frac{41}{12}$, $q(f_2)=-\frac{59}{24}$ modulo $2$ and scalar product $f_1.f_2=\frac{7}{4}$ modulo $1$.

In order to prove that the transcendental lattice corresponds to the Gram matrix $\begin{pmatrix}
	12 & 0\\
	0 & 24\\
\end{pmatrix}$
we must find for the corresponding quadratic form generators $g_1$ and $g_2$ satisfying $q(g_1)=\frac{41}{12}$, $q(g_2)=\frac{59}{24}$ modulo $2$ and scalar product $g_1.g_2=-\frac{7}{4}$ modulo $1$.
This is obtained with $g_1=\begin{pmatrix}
	\frac{1}{4}\\
	\frac{1}{3}\\
\end{pmatrix}$ and $g_2=\begin{pmatrix}
	\frac{5}{12}\\
	\frac{1}{8}\\
\end{pmatrix}$.
Thus $F_p$, hence $F_m$, are elliptic fibrations of the Kummer surface $K_{10}$.

\item
Self-isogenies of elliptic fibrations are quite interesting since they allow us to determine the Mordell-Weil lattice of the fibration. We are going to show two examples of that property. In both examples we know by specialisation an infinite section of the Mordell-Weil lattice of the fibration. The self-isogeny allows us to obtain another infinite section and achieve the determination of the Mordell-Weil lattice.

\begin{enumerate}
\item  From the elliptic fibration $F_{15}$ of $Y_{10}$
\[ (F_{\#15}) \qquad y^2=x(x+t^3-24t^2)(x+t^3-20t^2+4t),\]
with the elliptic parameters $\frac{x}{t^2(t-24)}$, $\frac{x}{t^2(t+1)}$, $\frac{x}{t_1(t_1-4)(t_1-24)}$, we get successively
\[(E_2) \qquad y^2=x^3-4t(t+1)(6t+5)x^2+4t^2(t+1)^3x,\]
\[(E_3) \qquad y^2=x^3-2t(t^2-14t-2)x^2+t^4(t-4)(t-24)x,\]
\[(\tilde{E_4}) \qquad y^2=x^3-28t_1^2(t_1-1)x^2+4t_1^3(t_1-1)^2(24t_1+1)x,\]

then letting $t_1=\frac{2t}{t+1}$ we get the fibration $E_4$ of $Y_{10}$

\begin{align*}
  (E_4) \qquad Y^{2}  & =X^{3}-112t^{2}\left(  t^{2}-1\right)  X^{2}+32t^{3}\left(
t^{2}-1\right)  ^{2}\left(  49t+1\right)  X\\
& III^{\ast}\left(  0\right)  ,\quad2I_{0}^{\ast}\left(  \pm1\right)  \quad
I_{2}\left(  \frac{-1}{49}\right)  \quad I_{1}\left(  \frac{1}{49}\right).
\end{align*}
We can see easily that the $2$-isogenous elliptic curve $EE_4$ of $E_4$ given by the equation
\[(EE_4) \qquad Y^{2}   =X^{3}+224t^{2}\left(  t^{2}-1\right)  X^{2}+128t^{3}\left(
    t^{2}-1\right)  ^{2}\left(  49t-1\right)  X\]
is isomorphic to $E_4$, the isomorphism from $EE_4$ to $E_4$ given by $t\mapsto-t,$ $X\mapsto 2X,$
and $Y \mapsto i2\sqrt{2}Y.$

The infinite section $P=(192t^{3}\left(  t+1\right)  ,32\sqrt{6}\left(  t+1\right)
^{2}t^{3}\left(  1+23t\right)  )$ on $E_{4} $ gives by $2$-isogeny  \ref{2-iso} an infinite section $P_1$ on $EE_4$ defined as
\begin{align*}
P & =\left(  X(t),Y(t)\right)  \mapsto P_{1}=\left(  X_{1}\left(  t\right)
,Y_{1}\left(  t\right)  \right)  \\
X_{1}\left(  t\right)    & =\frac{Y\left(  t\right)  ^{2}%
}{X\left(  t\right)  ^{2}},\quad Y_{1}\left(  t\right)  =-\frac{Y\left(  t\right)  \left(  X\left(  t\right)  ^{2}-B\right)
}{X\left(  t\right)  ^{2}}\\
B  & =32t^{3}\left(  t^{2}-1\right)  ^{2}\left(  49t+1\right).
\end{align*}
Finally by the previous isomorphism we get an infinite section $Q$ on $E_4$
\[
Q=(\frac{-1}{12}\left(  t-1\right)  ^{2}\left(  23t-1\right)  ^{2}%
,\frac{i\sqrt{3}}{72}\left(  t-1\right)  ^{2}\left(  23t-1\right)  \left(
1103t^{3}-97t^{2}-47t+1\right)).
\]
Computing the height matrix of $P$ and $Q$ we see that these infinite sections generate the
Mordell-Weil lattice of the elliptic fibration $E_4$ of $Y_{10}$.
\item Another example is obtained from $E_9$,
  \begin{align*}
E_{9}  & :Y^{2}=X^{3}+28t^{2}X^{2}+t^{3}\left(  t^{2}+98t+1\right)  X\\
& 2III^{\ast}\left(  0,\infty\right)  ,\quad2I_{2}\left(  t^{2}+98t+1\right)
,\quad2I_{1}\left(  t^{2}-98t+1\right).
\end{align*}

The isomorphism is given by $t\mapsto-t,$ $X\mapsto-\frac{1}{2}X,$ and
$Y\mapsto-\frac{1}{i2\sqrt{2}}Y.$

The specialized section is
\[
P_{9}=\left(  \frac{1}{96}\left(  t-1\right)  ^{2}\left(  t^{2}+98t+1\right)
,\frac{\sqrt{6}}{2304}\left(  t^{2}-1\right)  \left(  t^{2}+98t+1\right)
\left(  t^{2}+46t+1\right)  \right)
\]
and the image by complex multiplication is
\begin{align*}
Q_{9}  & =\left(  X_{9},Y_{9}\right)  \\
X_{9}  & = -\frac{1}{192}\frac{\left(  t-1\right)  ^{2}\left(  t^{2}%
-46t+1\right)  ^{2}}{\left(  t+1\right)  ^{2}}\\
Y_{9}  & =-\frac{i\sqrt{3}}{4608}\frac{\left(  t-1\right)  \left(
t^{2}-46t+1\right)  \left(  t^{6}-94t^{5}-385t^{4}+8636t^{3}-385t^{2}%
-94t+1\right)  }{\left(  t+1\right)  ^{3}}%
\end{align*}

Computing the height matrix of $P_9$ and $Q_{9}$ we get the matrix $[3 \quad 0 \quad 6]$ hence  $P_9$ and $Q_9$ generate the
Mordell-Weil lattice of the fibration $E_9$ of $Y_{10}$.

\end{enumerate}

\end{enumerate}
\end{proof}
\begin{remark}

  The $K3$ surfaces $Y_2$ and $Y_{10}$ are CM by $\mathbb Q(\sqrt{-2})$ and belong by their fibration $\#9$ to the family exhibited by van Geemen and Sch\"{u}tt in \cite{GS} Proposition 6.2.
  Indeed fibration $\#9$ of $Y_2$ is given in Bertin and Lecacheux paper \cite{BL2} as
  \[y^2=x^3+4g^2x^2+g^3(g+1)^2x\]
  with $\alpha(t^2)=2t^2$ and $\beta(t)=t^2+t$.

  Fibration $\#9$ of $Y_{10}$ is given in the present paper as $E_9$, thus with  $\alpha(t^2)=14t^2$ and $\beta(t)=t^2+t$. 
\end{remark}

\section{Mordell-Weil lattice of the specialised \#16 fibration of $Y_{10}$}

The rank of the specialisation for $k=10$ of elliptic fibrations increases
by one, so we have to determine one more generator for the Mordell-Weil
group. We give an example where the computation is easy using a two-isogeny
between an elliptic fibration of $Y_{10}$ and an elliptic fibration of the
Kummer surface $K_{10}=\text{Kum}(E_1,E_2)$ associated to $Y_{10}$, where $%
E_1,E_2$ are elliptic curves with complex multiplication. Then using the
method developped in \cite{Sh2} and \cite{KU}, we determine a section
on an elliptic fibration
of $Y_{10}.$ From \cite{BL2} Corollary 4.1, the two elliptic curves $E_{1}$
and $E_{2}$ have respective invariants $j_{1}=8000$ and $j_{2}=188837384000-
77092288000\sqrt{6}.$ Take 
\begin{equation*}
E_{1}:Y^{2}=X\left( X^{2}+4X+2\right)
\end{equation*}
as a model of the first curve. The $2$-torsion sections have $X$-coordinates 
$0,$ $-2\pm\sqrt{2}$.
The elliptic curve $E_{1}$ has complex
multiplication by $m_{2}=\sqrt{-2}$ defined by 
\begin{equation*}
\left( X,Y\right) \overset{m_{2}}{\mapsto}\left( -\frac{1}{2}\frac
{X^{2}+4X+2}{X},\frac{i\sqrt{2}}{4}\frac{Y\left( X^{2}-2\right) }{X^{2}}%
\right) .
\end{equation*}
Let $C_{3}$ and $\widetilde{C_{3}}$ the two groups of order $3$ generated by
the points of respective $X$-coordinates $\frac{1}{3}\left( 1+i\sqrt
{2}\right) $ and $\frac{1}{3}\left( 1-i\sqrt{2}\right) .$ These groups are
fixed by $m_{2}$ while the two order $3$ groups $\Gamma_{3}$ and $\widetilde{%
\Gamma_{3}}$ generated by the points of respective $X$-coordinates $-1+\sqrt
{6}$ and $-1-\sqrt{6}$ are exchanged by $m_{2}.$

If $M=\left( X_M,Y_M\right) $ is a generic point of $E_{1}$, we determine the image of $M$ by the isogeny of kernel $\Gamma_{3}$ with $%
X_{2}=\sum_{S\in \Gamma_3}X_{M+S}+k$ and $Y_{2}=\sum_{S\in \Gamma_3}Y_{M+S}$
where $k$ can be chosen so that the image of $\left( 0,0\right) $ is $%
X_{2}=0.$ It follows the $3$-isogeny

\begin{eqnarray*}
w_{3}:X_{2}&=&\phi _{3}\left( X\right) =\frac{X\left( X-2-\sqrt{6}\right) ^{2}%
}{\left( X+2-\sqrt{6}\right) ^{2}}, \\ 
Y_{2}&=&Y\psi _{3}\left( X\right) =-%
\frac{Y\left( X^{2}+\left( 8-2\sqrt{6}\right) X+2\right) \left( X-2-\sqrt{6}%
\right) }{\left( X+2-\sqrt{6}\right) ^{3}}
\end{eqnarray*}
and its $3$-isogenous curve $E_{2}$

\begin{align*}
E_{2}& :Y_{2}^{2}=X_{2}^{3}+28X_{2}^{2}+\left( 98+40\sqrt{6}\right) X_{2} \\
j\left( E_{2}\right) & =188837384000-77092288000\sqrt{6}.
\end{align*}

We use also $w_{6}=w_{3}\circ m_{2}$%
\begin{eqnarray*}
X_{2} &=&\phi _{6}\left( X\right) =-\frac{1}{2}\frac{\left(
X^{2}+4X+2\right) \left( X^{2}+2X\left( 4+\sqrt{6}\right) +2\right) ^{2}}{%
X\left( X+2+\sqrt{6}\right) ^{2}\left( X-2+\sqrt{6}\right) ^{2}} \\
Y_{2} &=&Y\psi _{6}\left( X\right) =\\
&&\frac{i\sqrt{2}}{4}\frac{Y\left( X^{2}-2\right) \left( X^{2}+2X\left( 4+%
\sqrt{6}\right) +2\right) H_4(X) 
}{X^{2}\left( X+2+\sqrt{6}\right) ^{3}\left( X-2+\sqrt{6}\right) ^{3}} \\
H_4&=&\left( X^{4}+4\left( \sqrt{6}-2\right)
X^{3}+4\left( -9+\sqrt{6}\right) X^{2}+8\left( \sqrt{6}-2\right) X+4\right).
\end{eqnarray*}

An equation for the Kummer surface $K_{10}$ is therefore 
\begin{equation*}
K_{10}:X\left( X^{2}+4X+2\right) =y^{2}X_{2}\left( X_{2}^{2}+28X_{2}+98+40%
\sqrt{6}\right) .
\end{equation*}

\subsubsection{Elliptic fibrations of $K_{10}$ and $Y_{10}$}

We use the following units of $\mathbb{Q}\left( \sqrt{2},\sqrt{3}\right) $%
\begin{align*}
r_{1} & =1+\sqrt{2}+\sqrt{6},\quad r_{1}^{\prime}=1-\sqrt{2}+\sqrt{6},\quad
        r_{1}r_{1}^{\prime}=s=\left( \sqrt{2}+\sqrt{3}\right) ^{2} \\
  r_{2} & =1-2\sqrt{3}-\sqrt{6},\quad r_{2}^{\prime}=1+2\sqrt{3}-\sqrt{6}
\quad r_{2}r_{2}^{\prime}=-s.
\end{align*}

In this paragraph we construct an elliptic fibration of $K_{10}$ giving after a two-isogeny the specialisation of the elliptic fibration $\#16$ on $Y_{10}$.

We consider the fibration 
\begin{align*}
K_{10} & \rightarrow\mathbb{P}^{1} \\
\left( X,X_{2},y\right) & \mapsto t=\frac{X_{2}}{X}.
\end{align*}

Notice that $X_{2}=tX$ and a Weierstrass equation $K_t$ ($K_t$ denotes $%
(K_{10})_t$ for simplification) for this fibration is obtained with the
following transformation

\begin{equation*}
X=-\sqrt{2}\left( 1+\sqrt{2}\right) \frac{X_{1}-2t\left( t-r_{1}^{2}\right)
\left( t-r_{2}^{2}\right) }{X_{1}\left( 3+2\sqrt{2}\right) -2t\left(
t-r_{1}^{2}\right) \left( t-r_{2}^{2}\right) },\qquad y=2\sqrt{2}\frac{X_{1}%
}{Y_{1}}
\end{equation*}

\begin{eqnarray}
K_{t} &:&Y_{1}^{2}=X_{1}\left( X_{1}-2t\left( t-r_{1}^{2}\right) \left(
t-r_{1}^{\prime 2}\right) \right) \left( X_{1}-2t\left( t-r_{2}^{2}\right)
\left( t-r_{2}^{\prime 2}\right) \right)  \label{A} \\
&=&f_{t}\left( X_{1}\right) .
\end{eqnarray}

The singular fibers are in $t=0$ and $\infty$ of type $I_{2}^{\ast}$ and at $%
t=r_{1}^{2},r_{2}^{2},r_{1}^{\prime2},r_{2}^{\prime2}$ of type $I_{2}.$ The
rank of the Mordell-Weil group is $2.$

\subsubsection{Sections on the elliptic fibration $K_t$}


We find two sections on the previous fibration using $w_{3}$ and $w_{6}\in
\hom \left( E_{1},E_{2}\right) .$ We consider the graph of $w_{3}$ and $w_{6}
$ on $E_{1}\times E_{2}$ and its image on $K_{10}=E_{1}\times E_{2}/\pm 1$ ,
then we have the multisections $Q_{3}$ and $Q_{6}$ $\ $\ defined
respectively by $\left( t=\frac{\phi _{3}\left( X\right) }{X},y=\frac{1}{%
\psi _{3}\left( X\right) }\right) $ and $\left( t=\frac{\phi _{6}\left(
X\right) }{X},y=\frac{1}{\psi _{6}\left( X\right) }\right) .$ In the
Weierstrass form this gives $Q_{3}=\left( X_{1,3},Y_{1,3}\right) $ where $%
X_{1,3}$ is a root of an irreducible polynomial $g_{3}\left( X_{1}\right)
\in \mathbf{C}\left( t\right) \left[ X_{1}\right] $ and $Y_{1,3}=h_{3}\left(
X_{1}\right) \in \mathbf{C}\left( t\right) \left[ X_{1}\right] \mathbf{.}$ \
In the same way we compute $g_{6}$ and $h_{6}$ for $Q_{6}.$ To have a
section on the fibration of parameter $t$ we compute the trace $P_{3}$ of $%
Q_{3}$ and $P_{6}$ of $Q_{6}.$ Recall for $Q_{3}$ the algorithm given in
\cite{BE} Algorithm $1:$ 

$g=g_{3}$, $h=h_{3}\mod g$. While $\deg \left( g\right) >1$ do $%
g:=\left( h^{2}-f_{t}\right) /g$,\quad\ $h:=h$ mod $g$ endwhile return $%
\left( g,h\right) $

The final results are 
\begin{align*}
P_{3}& =\left( x_{P_{3}},y_{P_{3}}\right)  \\
x_{P_{3}}& =\frac{1}{2s}\left( t+s\right) ^{2}\left( t-r_{1}^{2}\right)
\left( t-r_{1}^{\prime 2}\right) ,\quad y_{P_{3}}=x_{P_{3}}\frac{\left( 
\sqrt{6}-2\right) }{2}\left( \frac{t-s}{t+s}\right) \left( t^{2}-14t-4\sqrt{6%
}t+s^{2}\right) 
\end{align*}

and $P_{6}$ 
\begin{eqnarray*}
P_{6} &=&\left( x_{P_{6},}y_{P_{6}}\right) \\
x_{P_{6}} &=&\frac{-8\left( 5+\sqrt{6}\right) ^{2}\left( 19t^{2}-\left(
326+140\sqrt{6}\right) t+19s^{2}\right) ^{2}}{361\left( t-1\right)
^{2}\left( t-s^{2}\right) ^{2}} \\
y_{P_{6}} &=&\frac{i4\sqrt{2}\left( 5+\sqrt{6}\right) t^{2}\left(
t^{2}-s^{2}\right) \left( 19t^{2}-\left( 326+140\sqrt{6}\right)
t+19s^{2}\right) \left( \left( t^{2}-2t+s^{2}\right) ^{2}-96t^{2}\right) }{%
19\left( t-1\right) ^{3}\left( t-s^{2}\right) ^{3}}
\end{eqnarray*}

two sections on the fibration of the Kummer surface $K_{10}$. \ \ 

\subsubsection{Sections on the fibration \#16 of $Y_{10}$}

The $2$-isogenous elliptic curve of (\ref{A}) in the isogeny of kernel $%
\left( 0,0\right) $ has a Weierstrass equation 
\begin{align}
Y_{3}^{2} & =X_{3}\left( X_{3}^{2}+8\,t\left( {t}^{2}-28\,t+{s}^{2}\right)
X_{3}+64\,{\frac{{t}^{4}}{{s}^{2}}}\right)  \label{B} \\
X_{3} & =\left( \frac{Y_{1}}{X_{1}}\right) ^{2},\qquad Y_{3}=\frac
{Y_{1}(B-X_{1}^{2})}{{X_{1}}}  \notag
\end{align}
where $B$ is the coefficient of $X_{1}$ in (\ref{A}). Singular fibers are in 
$t=0$ and $\infty$ of type $I_{4}^{\ast}$ and of type $I_{1}$ at $%
t=r_{1}^{2},r_{2}^{2},r_{1}^{\prime2},r_{2}^{\prime2}$ $.$

We can show that it is a fibration of $Y_{10}$ . More precisely, the
specialisation of the fibration $\#16$ (\cite{BL2} Table 4) for $k=10$ has
a Weierstrass equation 
\begin{equation*}
y^{2}=x^{3}+t_{0}\left( 4\left( t_{0}^{2}+s^{2}\right) +t_{0}\left(
s^{4}+14s^{2}+1\right) x^{2}\right) +16s^{6}t_{0}^{4}x
\end{equation*}
with parameter $t_{0}$ . If $t=-t_{0}s^{2}$ and after scaling ($%
x=-s^{6}X_{3}/2,y=is^{9}Y_{3}/2^{3/2}$) we recover the equation (\ref{B}).
The image of $P_{3}$ by the isogeny, in the Weierstrass equation (\ref{B})\
is $p_{3}=\left( \xi _{p_{3}},\eta _{p_{3}}\right) $ with 
\begin{align*}
\xi _{p_{3}}& =\frac{1}{2s}\frac{\left( t^{2}-14t-4\sqrt{6}t+s^{2}\right)
^{2}\left( t-s\right) ^{2}}{(t+s)^{2}} \\
\eta _{p_{3}}& =-\frac{\left( -2+\sqrt{6}\right) }{4s}\frac{\left(
t^{2}-14t-4\sqrt{6}t+s^{2}\right) \left( t-s\right) L_{t}}{(t+s)^{3}} \\
\text{where}\quad L_{t}& =t^{6}+2\left( 1+2\sqrt{6}\right) t^{5}-\left(
993+404\sqrt{6}\right) t^{4}+\left( 17820+7272\sqrt{6}\right) t^{3} \\
& -\left( k97137+39656\sqrt{6}\right) t^{2}+\left( 56642+23124\sqrt{6}%
\right) t+s^{6}.
\end{align*}

The image of $P_{6}$ is $2-$divisible \ and equal to $2p_{6}$ with 
\begin{eqnarray*}
p_{6} &=&\left( \xi _{6},\eta _{6}\right)  \\
\xi _{6} &=&\frac{-8t^{3}\left( t-1\right) ^{2}}{\left( t-s^{2}\right) ^{2}}
\\
\eta _{6} &=&\frac{32}{19}\frac{i\sqrt{2}\left( 5+\sqrt{6}\right)
t^{4}\left( t-1\right) \left( 19t^{2}-\left( 326+140\sqrt{6}\right)
t+19s^{2}\right) }{\left( t-s^{2}\right) ^{3}}
\end{eqnarray*}

We verify using definitions that $<p_{3},p_{6}>=0$ and $ht\left(
p_{3}\right) .ht\left( p_{6}\right) =18$. So by Shioda-Tate formula (\cite
{Sh3}, Corollary 1.7) $p_3$, $p_6$ and $(0,0)$ generate the
Mordell-Weil group.

\section{$3$-isogenies from elliptic fibrations of $Y_2$}\label{sec:2}
 Let us prove the following theorem 
  \begin{theorem}\label{3isog}
Given a Weierstrass equation $E_t$ of a $K3$ surface $S$ with a $3$-torsion section given by the point $(0,0)$

\[(E_t) \qquad Y^2+A(t)XY+B(t)Y=X^3\]
and transcendental lattice $T_S$, denote $\tilde{S}$ the $3$-isogenous $K3$ surface $S/<\omega>$ where $\omega=(0,0)$ and $T_{\tilde{S}}$ its transcendental lattice.

A sufficient condition to get $T_{\tilde{S}}=T_S[3]$ is that $Y$ will be the elliptic parameter for an elliptic fibration with two $II^*$ fibers at $0$ and $\infty$. In particular this condition can be realized if $(E_t)$ has singular fibers either $I_{12}$and $IV^*$ or $I_{18}$.

\end{theorem}

\begin{proof}
  Consider first the case of $(E_t)$ with a singular fiber $I_{18}$. From the components $\theta_i$ we can exhibit two fibers of type $II^*$.

\newpage

  \begin{multicols}{2}
    \unitlength 1cm
\begin{center}
\begin{tikzpicture}[scale=0.8]

\thinlines 

\draw [fill=black](9,0) circle (0.07cm)node [above] {$\theta_{\infty,0}$};
\draw [fill=black](8,1) circle (0.07cm)node [above] {$\theta_{\infty,1}$};
\draw [fill=black](7,2) circle (0.07cm)node [above] {$\theta_{\infty,2}$};
\draw [fill=black](6,3) circle (0.07cm)node [above] {$\theta_{\infty,3}$};
\draw [fill=black](5,3) circle (0.07cm)node [above] {$\theta_{\infty,4}$};
\draw [fill=black](4,3) circle (0.07cm)node [above] {$\theta_{\infty,5}$};
\draw [fill=black](3,3) circle (0.07cm)node [above] {$\theta_{\infty,6}$};
\draw [fill=black](2,2) circle (0.07cm)node [above] {$\theta_{\infty,7}$};
\draw [fill=black](1,1) circle (0.07cm)node [above] {$\theta_{\infty,8}$};
\draw [fill=black](0,0) circle (0.07cm)node [above] {$\theta_{\infty,9}$};
\draw [fill=black](1,-1) circle (0.07cm)node [above] {$\theta_{\infty,10}$};
\draw [fill=black](2,-2) circle (0.07cm)node [above] {$\theta_{\infty,11}$};
\draw [fill=black](3,-3) circle (0.07cm)node [below] {$\theta_{\infty,12}$};
\draw [fill=black](4,-3) circle (0.07cm)node [below] {$\theta_{\infty,13}$};
\draw [fill=black](5,-3) circle (0.07cm)node [below] {$\theta_{\infty,14}$};
\draw [fill=black](6,-3) circle (0.07cm)node [below] {$\theta_{\infty,15}$};
\draw [fill=black](7,-2) circle (0.07cm)node [above] {$\theta_{\infty,16}$};
\draw [fill=black](8,-1) circle (0.07cm)node [above] {$\theta_{\infty,17}$};

\draw (9,0)--(6,3);
\draw (6,3)--(3,3);
\draw (3,3)--(0,0);
\draw (0,0)--(3,-3);
\draw (3,-3)--(6,-3);
\draw (6,-3)--(9,0);

\draw [fill=red](7,0) circle (0.07cm)node [above] {$O$};
\draw (9,0)--(7,0);
\draw [fill=red](4.1,1.5) circle (0.07cm)node [above] {$T_3$};
\draw  (4.1,1.5)--(3,3);
\draw [fill=red](4.1,-1.5) circle (0.07cm)node [above] {$2T_3$};
\draw(4.1,-1.5)--(3,-3);
\end{tikzpicture}%
\end{center}

\columnbreak

We consider the two divisors $D$
\begin{eqnarray*}
D&=&3(O)+6\theta _{\infty ,0}+5\theta _{\infty ,17}+4\theta _{\infty
     ,16}+3\theta _{\infty ,15}\\
  &&+2\theta _{\infty ,14}+\theta _{\infty ,13} +4\theta _{\infty ,1}+2\theta _{\infty ,2}
\end{eqnarray*}
and $D'$ 
\begin{eqnarray*}
D^{\prime }&=&3(T_{3})+6\theta _{\infty,6}+5\theta _{\infty,7}+4\theta
               _{\infty,8}+3\theta _{\infty ,9}\\
  &&+2\theta _{\infty ,10}+\theta _{\infty ,11} +4\theta _{\infty,5}+2\theta _{\infty,4}
\end{eqnarray*}
\end{multicols}

We have $D.D'=0$, $D.\theta_{\infty,12}=1$, so $D$ and $D'$ are two singular fibers of a new elliptic fibration 
with two fibers of type $II^*$. The parameter $h$ of this new fibration can be chosen as $h=Y$ 
since the horizontal divisor of $D'-D$ is $3(T_3)-3(0)$ (\cite{L} Propositions 8.1 and 8.2) .

Since $E_t$ has only a singular fiber $I_{18}$ at $t=\infty$, it follows that $B(t)=b$ and $A(t)$ is a degree two polynomial.  After some scaling we can take $A(t)=t^2+a$ and $b=1$.
With $Y=h$ as new parameter and putting $X=-\frac{h^2(h+1)}{x}$, $t=-\frac{y}{hx}$, it follows the Weierstrass equation
\[y^2=x^3-ah^2x^2+h^5(h+1)^2\]
with singular fibers $II^*(0,\infty)$, $I_2(-1)$, $6I_1$. Proceeding as in Shioda's theorem \ref{th:sh} with $h=u^3$ we get the fibration $E_u$ of $S^{(3)}$
with transcendental lattice $T_S[3]$,
\[(E_u)\qquad Y^2=X^3-au^2X^2+u^3(u^3+1)^2,\]\label{Eu}
which is a rank $7$ elliptic fibration of $S^{(3)}$ with singular fibers $2I_0^*(\infty,0)$, $3I_2(u^3+1)$, $6I_1$.

The same transformation can be performed directly on the Weierstrass equation $(E_t)$ taking $Y=u^3$ and giving, if $X=ux$, the following cubic equation of  $S^{(3)}$
\[u^3+(t^2+a)ux+1=x^3.\]
Finally applying 2.7.3,  it follows a Weierstrass equation $(F_t)$ of this cubic $3$-isogenous to $(E_t)$ i.e. $\tilde{S}=S^{(3)}$
\[(F_t) \quad y^2-3(t^2+a)xy+(27-(t^2+a)^3)y=x^3.\]
\begin{remark}
Using Shimada \cite{Shim2}, we notice that the elliptic curve $(F_t)$ with singular fibers $I_6$ and $6I_3$ has $(\mathbb Z/3\mathbb Z)^2)$-torsion.
  
 \end{remark} 
\newpage

Consider now the case of $(E_t)$ with singular fibers $I_{12}(\infty)$ and $IV^*(0)$ and draw two fibers of type $II^*$.
\begin{multicols}{2}
\unitlength 1cm
\begin{center}
\begin{tikzpicture}[scale=0.8]

\thinlines 

\draw [fill=black](9,0) circle (0.07cm)node [above] {$\theta_{\infty,0}$};
\draw [fill=black](7.5,1.5) circle (0.07cm)node [above] {$\theta_{\infty,1}$};
\draw [fill=black](6,3) circle (0.07cm)node [above] {$\theta_{\infty,2}$};
\draw [fill=black](4.5,3) circle (0.07cm)node [above] {$\theta_{\infty,3}$};
\draw [fill=black](3,3) circle (0.07cm)node [above] {$\theta_{\infty,4}$};
\draw [fill=black](1.5,1.5) circle (0.07cm)node [above] {$\theta_{\infty,5}$};
\draw [fill=black](0,0) circle (0.07cm)node [above] {$\theta_{\infty,
6}$};
\draw [fill=black](1.5,-1.5) circle (0.07cm)node [above] {$\theta_{\infty,7}$};
\draw [fill=black](3,-3) circle (0.07cm)node [below] {$\theta_{\infty,8}$};
\draw [fill=black](4.5,-3) circle (0.07cm)node [below] {$\theta_{\infty,9}$};
\draw [fill=black](6,-3) circle (0.07cm)node [below] {$\theta_{\infty,10}$};
\draw [fill=black](7.5,-1.5) circle (0.07cm)node [above] {$\theta_{\infty,11}$};

\draw (9,0)--(6,3);
\draw (6,3)--(3,3);
\draw (3,3)--(0,0);
\draw (0,0)--(3,-3);
\draw (3,-3)--(6,-3);
\draw (6,-3)--(9,0);

\draw [fill=black](6,0) circle (0.07cm)node [above] {$\theta_{0,4}$};
\draw [fill=red](8.2,0) circle (0.07cm)node [above] {$O$};
\draw (8.2,0) circle (0.11cm);
\draw [fill=black](7.5,0) circle (0.07cm)node [above] {$\theta_{0,0}$};
\draw [fill=black](6.7,0) circle (0.07cm)node [above] {$\theta_{0,1}$};
\draw (9,0)--(6,0);
\draw [fill=red](3.75,2.25) circle (0.07cm)node [above] {$T_3$};
\draw (3.75,2.25) circle (0.11cm);
\draw [fill=black](4.5,1.5) circle (0.07cm)node [above] {$\theta_{0,2}$};
\draw [fill=black](5.25,0.75) circle (0.07cm)node [above] {$\theta_{0,3}$};
\draw  (6,0)--(3,3);
\draw [fill=black](5.25,-0.75) circle (0.07cm)node [above] {$\theta_{0,5}$};
\draw [fill=black](4.5,-1.5) circle (0.07cm)node [above] {$\theta_{0,6}$};
\draw [fill=red](3.75,-2.25) circle (0.07cm)node [above] {$2T_3$};
\draw (3.75,-2.25) circle (0.11cm);
\draw(6,0)--(3,-3);
\end{tikzpicture}%
\end{center}

\columnbreak

Consider the divisors
$D$ \ 
\begin{eqnarray*}
D &=&3(O)+6\theta _{\infty ,0}+5\theta _{\infty ,11}+4\theta _{\infty
      ,10}+3\theta _{\infty ,9}\\
  &&+2\theta _{\infty ,8}+\theta _{\infty ,7} +4\theta _{\infty ,1}+2\theta _{\infty ,2}
\end{eqnarray*}
and $D^{\prime }$%
\begin{eqnarray*}
D^{\prime } &=&3\theta _{0,1}+6\theta _{0,4}+5\theta _{0,3}+4\theta
                _{0,2}+3(T_{3})\\
  &&+2\theta _{\infty ,4,}+\theta _{\infty ,5} +4\theta _{0,5}+2\theta _{0,6}
\end{eqnarray*}

\end{multicols}


They satisfy $D.D'=0$, $D.\theta_{\infty,6}=1$  and are two singular fibers of a new elliptic fibration 
with two fibers of type $II^*$. The parameter $h$ of this new fibration can be chosen as $h=Y$ 
since the horizontal divisor of $D'-D$ is $3(T_3)-3(0)$ (\cite{L} Propositions 8.1 and 8.2).

Moreover a $K3$ surface $S$ having an elliptic fibration with $3$-torsion and two fibers $I_{12}(\infty)$ and $IV^*(0)$ has a Weierstrass equation of the form
\[(E) \qquad y^2-t(t+l)xy+t^2y=x^3.\]
Putting $y=h$ we obtain a cubic in $t$ and $x$ with a point $(x=1,t=\frac{h^2-1}{hl})$; thus
 it is an elliptic curve with a Weierstrass form 
\[Y^2=X^3-(1/48l^4+3)t^4X+t^7-1/864l^2(-648+l^4)t^6+t^5.
\] 
With the base change $h=u^3$ we get a Weierstrass equation of an elliptic fibration of the $K3$ surface $S^{(3)}$ with transcendental lattice $T_{S^{(3)}}=T_S[3]$.

The same transformation can be performed directly on the Weierstrass equation $(E)$ taking $y=u^3$ and giving, if $x=uX$, the following equation of a cubic on $S^{(3)}$
\[u^3-t(t+l)uX+t^2=X^3.\]
Finally applying 2.7.3  it follows a Weierstrass equation of this cubic $3$-isogenous to $(E)$
\[Y^2+3t(t+l)YX+t^2(t^4+3t^3l+3t^2l^2+tl^3+27)Y=X^3.
\]

\end{proof}

\begin{theorem}\label{th:spe}
  \begin{enumerate}
  \item The $3$-isogenous elliptic fibrations $H_{\#19}(k)$ (resp.$H_{\#20}(k)$) of the two generic fibrations $E_{\#19}(k)$ (resp. $E_{\#20}(k)$) are elliptic fibrations of the same $K3$ surface $N_k$ with transcendental lattice $U(3)\oplus <4>$, $U$ being the hyperbolic lattice.
  \item The specialized $K3$ surface $N_2$ is the $K3$ surface $Y_{10}$.
    \item The specialized $K3$ surface $N_{10}$ is the $K3$ surface with transcendental lattice $[4 \quad 0 \quad 18]$ (Shimada notation \cite{SZ}).
\end{enumerate}
  \end{theorem}
  
  \begin{proof}
\begin{enumerate}   
\item
The $6$-torsion elliptic fibration $\#20$ has a Weierstrass equation
\[E_{\#20}(k) \qquad y^2-(t^2-tk+3)xy-(t^2-tk+1)y=x^3\]
with singular fibers $I_{12}(\infty)$, $2I_3(t^2-kt+1)$, $2I_2(0,k)$, $2I_1(t^2-kt+9)$ and $3$-torsion point $(0,0)$. Using \ref{formulae}, it follows the Weierstrass equation of its $3$-isogenous fibration
\[H_{\#20}(k)=E_{\#20}(k)/\langle (0,0) \rangle  \qquad Y^2+3(t^2-tk+3)XY+t^2(t^2-tk+9)(t-k)^2Y=X^3\]
with singular fibers $2I_6(0,k)$, $I_4(\infty)$, $2I_3(t^2-kt+9)$, $2I_1(t^2-kt+1)$. Thus it is a rank $0$ and $6$-torsion elliptic fibration of a $K3$-surface with Picard number $19$ and discriminant $\frac{6\times 6 \times 3 \times 3 \times 4}{6\times 6}=12\times 3$.

Now we shall compute the Gram matrix $NS(20)$ of the N\'eron-Severi lattice of the $K3$ surface with elliptic fibration $H_{\#20}(k)$ in order to deduce its discriminant form.

Applying Shioda's result \ref{shio}, we order the following elements as, $s_0$, $F$, $\theta_{0,i}$, $1\leq i \leq 4$, $s_3$, $\theta_{k,i}$, $1\leq i \leq 5$, $\theta_{\infty,i}$, $1\leq i \leq 3$, $\theta_{t_0,i}$, $1\leq i \leq 2$,  $\theta_{t_1,i}$, $1\leq i \leq 2$, where $s_0$ and $s_3$ denotes respectively the zero and $3$-torsion section, $F$ the generic section, $\theta_{k,i}$ the components of reducible singular fibers, $t_0$ and $t_1$ being roots of $t^2-kt+9$. We obtain

\[ NS(20)=\left(
\begin{smallmatrix}
        -2&1&0&0&0&0&0&0&0&0&0&0&0&0
&0&0&0&0&0\\
        1&0&0&0&0&0&1&0&0&0&0&0&0&0&0&0&0&0&0  \\
        0&0&-2&1&0&0&1&0&0&0&0&0&0&0&0&0&0&0&0 \\
        0&0&1&-2&1&0&0&0&0&0&0&0&0&0&0&0&0&0&0  \\
        0&0&0&1&-2&1&0&0&0&0&0&0&0&0&0&0&0&0&0\\
        0&0&0&0&1&-2&0&0&0&0&0&0&0&0&0&0&0&0&0\\
        0&1&1&0&0&0&-2&1&0&0&0&0&0&1&0&0&1&0&1 \\
        0&0&0&0&0&0&1&-2&1&0&0&0&0&0&0&0&0&0&0 \\
        0&0&0&0&0&0&0&1&-2&1&0&0&0&0&0&0&0&0&0 \\                                                                                                                                                                  
        0&0&0&0&0&0&0&0&1&-2&1&0&0&0&0&0&0&0&0\\                                                                                                                                                                   
        0&0&0&0&0&0&0&0&0&1&-2&1&0&0&0&0&0&0&0 \\                                                                                                                                                                  
        0&0&0&0&0&0&0&0&0&0&1&-2&0&0&0&0&0&0&0     \\                                                                                                                                                              
        0&0&0&0&0&0&0&0&0&0&0&0&-2&1&0&0&0&0&0   \\                                                                                                                                                                
        0&0&0&0&0&0&1&0&0&0&0&0&1&-2&1&0&0&0&0 \\                                                                                                                                                                  
        0&0&0&0&0&0&0&0&0&0&0&0&0&1&-2&0&0&0&0    \\                                                                                                                                                               
        0&0&0&0&0&0&0&0&0&0&0&0&0&0&0&-2&1&0&0 \\                                                                                                                                                                  
        0&0&0&0&0&0&1&0&0&0&0&0&0&0&0&1&-2&0&0  \\                                                                                                                                                                 
        0&0&0&0&0&0&0&0&0&0&0&0&0&0&0&0&0&-2&1 \\                                                                                                                                                                  
        0&0&0&0&0&0&1&0&0&0&0&0&0&0&0&0&0&1&-2 \\                                                                                                                                                                  
\end{smallmatrix}                                                                                                                                                                                                  
\right ).\]

We get $\det(NS(20))=12\times 3$ and applying Shimada's lemma \ref{lem:gram}, the discriminant form $G_{NS(20)}\simeq \mathbb Z/3 \oplus \mathbb Z/12$ is generated by vectors $L_1$ and $L_2$ satisfying $q_{L_1}=0$, $q_{L_2}=-\frac{11}{12}$ and $b(L_1,L_2)=\frac{1}{3}$.
Denoting $M(20)$ the following Gram matrix of the lattice $U(3)\oplus \langle 4 \rangle$,
\[M(20)=\begin{pmatrix}                                                                                                                                                                                            
        0 & 0 & 3\\                                                                                                                                                                                                
        0 & 4 & 0\\                                                                                                                                                                                                
        3 & 0 & 0\\                                                                                                                                                                                                
\end{pmatrix},\]
we find for generators of its discriminant form the vectors
\[g_1=\begin{pmatrix}                                                                                                                                                                                              
        0\\                                                                                                                                                                                                        
        0\\                                                                                                                                                                                                        
        \frac{1}{3}\\                                                                                                                                                                                              
\end{pmatrix}              \qquad            g_2=\begin{pmatrix}                                                                                                                                                   
        \frac{1}{3}\\                                                                                                                                                                                              
        \frac{1}{4}\\                                                                                                                                                                                              
        \frac{1}{3}\\                                                                                                                                                                                              
\end{pmatrix}                                                                                                                                                                                                      
\]
satisfying $q_{g_1}=0$, $q_{g_2}=\frac{11}{12}$ and $b(g_1,g_2)=\frac{1}{3}$. We deduce that $M(20)$ is the transcendental lattice of the $K3$ surface with elliptic fibration $H_{\#20}(k)$.

A Weierstrass equation of the $3$-torsion, rank $1$, elliptic fibration $\#19$ can be written as
\[E_{\#19}(k) \qquad y^2+ktxy+t^2(t^2+kt+1)y=x^3\]
with singular fibers $2IV^*(0,\infty)$, $2I_3(t^2+kt+1)$, $2I_1(k^3t-27kt-27t^2-27)$ and infinite point $P=(-t^2,-t^2)$ of height $h(P)=\frac{4}{3}$.
Its $3$-isogenous elliptic fibration $E_{\#19}(k)/\langle (0,0) \rangle $ has a Weierstrass equation
\[H_{\#19}(k) \qquad Y^2-3ktXY-Yt^2(27t^2-k(k^2-27)t+27)=X^3.\]
It is a $3$-torsion, rank $1$, elliptic fibration of a $K3$-surface with Picard number $19$ and singular fibers $2IV^*(0,\infty)$, $2I_3(27t^2-k(k^2-27)t+27)$, $2I_1(t^2+kt+1)$ and infinite order point $Q$ with\
 $x$-coordinate $x_Q=-3-3kt-(k^2+3)t^2-3kt^3-3t^4$ and height $h(Q)=4$. This point $Q$ is the image of the point $P$ in the $3$-isogeny and non $3$-divisible, hence generator of the non torsion part of the Mordell-Weill lattice. We deduce the discriminant of this $K3$-surface $\frac{3\times 3 \times 3 \times 3\times 4}{3 \times 3}=12\times 3$.

Applying Shioda's result \ref{shio}, we order the components of the singular fibers as, $s_0$, $F$, $\theta_{0,i}$, $1\leq i \leq 6$, $\theta_{\infty,i}$, $1\leq i \leq 6$, $\theta_{t_0,i}$, $1\leq i \leq 2$, $s_3$, $\theta_{t_1,2}$, $s_{\infty}$,where $s_0$, $s_3$, $s_{\infty}$ denotes respectively the zero, $3$-torsion and infinite section, $F$ the generic section, $t_0$ and $t_1$ being roots of $27t^2-k(k^2-27)t+27$. The numbering of components of $IV^*$ is done using Bourbaki's notations \cite{Bour}. It follows the Gram matrix $NS(19)$ of the corresponding $K3$ surface

\[ NS(19)=                                                                                                                                                                                                         
\left (                                                                                                                                                                                                            
\begin{smallmatrix}                                                                                                                                                                                                
-2&1&0&0&0&0&0&0&0&0&0&0&0&0                                                                                                                                                                                       
&0&0&0&0&0\\1&0&0&0&0&0&0&0&0&0&0&0&0&0&0&0&1&0&1                                                                                                                                                                  
\\0&0&-2&0&0&1&0&0&0&0&0&0&0&0&0&0&0&0&0                                                                                                                                                                           
\\0&0&0&-2&1&0&0&0&0&0&0&0&0&0&0&0&1&0&0                                                                                                                                                                           
\\0&0&0&1&-2&1&0&0&0&0&0&0&0&0&0&0&0&0&0                                                                                                                                                                           
\\0&0&1&0&1&-2&1&0&0&0&0&0&0&0&0&0&0&0&0                                                                                                                                                                           
\\0&0&0&0&0&1&-2&1&0&0&0&0&0&0&0&0&0&0&0                                                                                                                                                                           
\\0&0&0&0&0&0&1&-2&0&0&0&0&0&0&0&0&0&0&0                                                                                                                                                                           
\\0&0&0&0&0&0&0&0&-2&0&0&1&0&0&0&0&0&0&0                                                                                                                                                                           
\\0&0&0&0&0&0&0&0&0&-2&1&0&0&0&0&0&1&0&0                                                                                                                                                                           
\\0&0&0&0&0&0&0&0&0&1&-2&1&0&0&0&0&0&0&0                                                                                                                                                                           
\\0&0&0&0&0&0&0&0&1&0&1&-2&1&0&0&0&0&0&0                                                                                                                                                                           
\\0&0&0&0&0&0&0&0&0&0&0&1&-2&1&0&0&0&0&0                                                                                                                                                                           
\\0&0&0&0&0&0&0&0&0&0&0&0&1&-2&0&0&0&0&0                                                                                                                                                                           
\\0&0&0&0&0&0&0&0&0&0&0&0&0&0&-2&1&0&0&0                                                                                                                                                                           
\\0&0&0&0&0&0&0&0&0&0&0&0&0&0&1&-2&1&0&0                                                                                                                                                                           
\\0&1&0&1&0&0&0&0&0&1&0&0&0&0&0&1&-2&1&2                                                                                                                                                                           
\\0&0&0&0&0&0&0&0&0&0&0&0&0&0&0&0&1&-2&0                                                                                                                                                                           
\\0&1&0&0&0&0&0&0&0&0&0&0&0&0&0&0&2&0&-2\                                                                                                                                                                          
\end{smallmatrix}                                                                                                                                                                                                  
\right ) .                                                                                                                                                                                                          
\]
Its determinant satisfies
$\det(NS(19))=12\times 3$ and according to Shimada's lemma \ref{lem:gram}, $G_{NS(19)}\simeq \mathbb Z/3 \oplus \mathbb Z/12$ is generated by vectors $M_1$ and $M_2$ satisfying $q_{M_1}=-\frac{2}{3}$, $q_{M_2}=\frac{5}{12}$ and $b(M_1,M_2)=0$. We find also generators for the transcendental discriminant form $M(20)$
\[h_1=\begin{pmatrix}                                                                                                                                                                                              
        \frac{1}{3}\\                                                                                                                                                                                              
        0\\                                                                                                                                                                                                        
        \frac{1}{3}\                                                                                                                                                                                               
\end{pmatrix}              \qquad            h_2=\begin{pmatrix}                                                                                                                                                   
        \frac{1}{3}\\                                                                                                                                                                                              
        \frac{1}{4}\\                                                                                                                                                                                              
        -\frac{1}{3}\\

\end{pmatrix}                                                                                                                                                                                                      
\]
satisfying $q_{h_1}=\frac{2}{3}$, $q_{h_2}=-\frac{5}{12}$ and $b(h_1,h_2)=0$. We deduce that $M(20)$ is also the transcendental lattice of the $K3$ surface with elliptic fibration $H_{\#19}(k)$.

It follows the discriminant form of its N\'eron-Severi lattice, 

$G_{NS(19)}=\mathbb Z/3\mathbb Z(-\frac{2}{3})\oplus \mathbb Z/12\mathbb Z(\frac{5}{12})$,
 which is also $G_{NS(20)}$ 
 since the generators $L'_1=15L_1+4L_2$, 
 $L'_2=L_1-L_2$ 
 satisfy $q_{L'_1}=-\frac{2}{3}$, 
 $q_{L'_2}=\frac{5}{12}$, $b(L'_1,L'_2)=0$.

\item

To prove that $N_2=Y_{10}$ it is sufficient to prove that $H_{\#20}(2)$ is an elliptic fibration of $Y_{10}$ since, by the previous theorem, $H_{\#19}(2)$ is another fibration of the same $K3$-surface.

But we see easily that $H_{\#20}(2)$ is the $6$-torsion extremal of $Y_{10}$ numbered $8$ in Shimada-Zhang \cite{SZ}.

\item Similarly, $N_{10}$ is the $K3$ surface with transcendental lattice $[4 \quad 0 \quad 18]$ since $H_{\#20}(10)$ is a fibration of that surface.
A Weierstrass equation for $H_{\#20}\left(  10\right)  $ is
\begin{align*}
&  Y^{2}+3\left(  t^{2}-22\right)  YX+\left(  t^{2}-25\right)  ^{2}\left(
t^{2}-16\right)  Y=X^{3}
\end{align*}
with singular fibers $I_{4}\left(  \infty\right)  ,2I_{6}\left(  \pm5\right)
,2I_{3}\left(  \pm4\right)  ,2I_{1}\left(  t^{2}-24\right)  .$ We have a
$2$-torsion section $s_2=\left(  -\left(  t^{2}-25\right)  ^{2},\left(
t^{2}-25\right)  ^{3}\right)  $ and a $6$-torsion section 

\noindent
$s_{6}=\left(  -\left(
t^{2}-25\right)  \left(  t^{2}-16\right)  ,-\left(  t^{2}-25\right)  \left(
t^{2}-16\right)  ^{2}\right)  $. 

\noindent
The section
$P_{t}=\left(  4\left(  t+5\right)  ^{2}\left(  t-4\right)  ,-\left(
t+5\right)  ^{4}\left(  t-4\right)  ^{2}\right)  $ is of infinite order and is
a generator of the Mordell-Weil lattice. \ Moreover we have $\theta_{\pm
5,1}.s_{6}=1,\theta_{\infty,2}.s_{6}=1,\theta_{\pm3,1}.s_{6}=1$. So the
following divisor is $6$-divisible%
\[
\sum_{i=1}^{5}\left(  6-i\right)  \theta_{\pm5,i}+3\theta_{\infty,1}%
+3\theta_{\infty,3}+4\theta_{\pm,4,1}+2\theta_{\pm4,2}\approx 6s_{6}%
\]
and we can replace $\theta_{5,5}$ by $s_{6}.$

Moreover we get $s_{6}.P_{t}=1$ (for $t=-3)$,
$\theta_{5,0}.P_{t}=1,\theta_{-5,4}.P_{t}=1$ and $\theta_{\infty,0}%
.P_{t}=1,\theta_{4,1}.P_{t}=1,\theta_{-4,0}.P_{t}=1.$ All these computations
give the Gram matrix of the N\'{e}ron-Severi lattice of discriminant $-72$

\[
\left (\begin {smallmatrix}
-2&1&0&0&0&0&0&0&0&0&0&0&0&0&0&0&0&0&0&0\\
1&0&0&0&0&0&0&0&0&0&0&1&0&0&0&0&0&0&0
&1\\0&0&-2&1&0&0&0&0&0&0&0&1&0&0&0&0&0&0&0&0
\\0&0&1&-2&1&0&0&0&0&0&0&0&0&0&0&0&0&0&0&0
\\0&0&0&1&-2&1&0&0&0&0&0&0&0&0&0&0&0&0&0&0
\\0&0&0&0&1&-2&1&0&0&0&0&0&0&0&0&0&0&0&0&0
\\0&0&0&0&0&1&-2&0&0&0&0&0&0&0&0&0&0&0&0&0
\\0&0&0&0&0&0&0&-2&1&0&0&1&0&0&0&0&0&0&0&0
\\0&0&0&0&0&0&0&1&-2&1&0&0&0&0&0&0&0&0&0&0
\\0&0&0&0&0&0&0&0&1&-2&1&0&0&0&0&0&0&0&0&0
\\0&0&0&0&0&0&0&0&0&1&-2&0&0&0&0&0&0&0&0&1
\\0&1&1&0&0&0&0&1&0&0&0&-2&0&1&0&1&0&1&0&1
\\0&0&0&0&0&0&0&0&0&0&0&0&-2&1&0&0&0&0&0&0
\\0&0&0&0&0&0&0&0&0&0&0&1&1&-2&1&0&0&0&0&0
\\0&0&0&0&0&0&0&0&0&0&0&0&0&1&-2&0&0&0&0&0
\\0&0&0&0&0&0&0&0&0&0&0&1&0&0&0&-2&1&0&0&1
\\0&0&0&0&0&0&0&0&0&0&0&0&0&0&0&1&-2&0&0&0
\\0&0&0&0&0&0&0&0&0&0&0&1&0&0&0&0&0&-2&1&0
\\0&0&0&0&0&0&0&0&0&0&0&0&0&0&0&0&0&1&-2&0
\\0&1&0&0&0&0&0&0&0&0&1&1&0&0&0&1&0&0&0&-2
\end {smallmatrix}\right ).
\]
According to Shimada's lemma \ref{lem:gram}, $G_{NS}=\mathbb{Z}/2\mathbb{Z}%
\oplus\mathbb{Z}/36$ is generated by the vectors $L_1$ and $L_2$ satisfying 
$q_{L_1}=\frac{-1}{2},q_{L_2}=\frac{-35}{36}$ and $b(L_1,L_2)=\frac{-1}{2}.$

Moreover the following generators of the discriminant group of the lattice with Gram matrix
$M_{18}=\left(\begin{smallmatrix}
4&0\\
0&18
\end{smallmatrix} \right)$, namely $f_1=(0,\frac{1}{2})$,$f_2=(\frac{1}{4},\frac{7}{18})$ verify $q_{f_1}=\frac{1}{2}$, $q_{f_2}=\frac{35}{36}$ and 
$b(f_1,f_2)=\frac{1}{2}.$
So the transcendental lattice is $[4 \quad 0 \quad 18]$.

\end{enumerate}
\end{proof}

\begin{theorem}
Define  $Y_{k}^{\left(  3\right)  }$ the elliptic surface  obtained by the
base change $\tau$ of the elliptic fibration of $Y_{k}$ with two singular
fibers of type $II^{\ast},$ where $\tau$ is the morphism given by $u\mapsto
h=u^{3}.$ Then 
the $K3$ surface $Y_{k}^{\left(  3\right)  }$ with transcendental lattice $U(3)\oplus<36>$ has a genus one fibration
without section such that its Jacobian variety 
 $J\left(  Y_{k}^{\left(  3\right)  }\right)  =N_{k}$ is the $K3$ surface with transcendental lattice $U(3)\oplus<4>$.
\end{theorem}
\begin{proof}
Recall a  Weierstrass equation for fibration $\#19$ of $Y_k$
\[
Y^{2}+tkYX+t^{2}\left(  t+s\right)  \left(  t+1/s\right)  Y=X^{3}%
\]
where $k=s+\frac{1}{s}.$ The fibration of $Y_{k}$ with two singular fibers
$II^{\ast}$ can be  obtained  with the parameter $h_{k}=\frac{Y}{\left(
t+s\right)  ^{2}}$ (\cite{BL2} Table 3). The surface  $Y_{k}^{\left(  3\right)  }$ is
  defined by $h=u^{3}$ and has then the following equation
\[
u^{3}s\left(  t+s\right)  +tkuWs+t^{2}\left(  ts-1\right)  -W^{3}s=0,
\]
where $X=\left(  t+s\right)  uW$.

We consider the fibration
\begin{align*}
Y_{k}^{\left(  3\right)  } &  \rightarrow\mathbb{P}^{1}\\
\left(  u,t,W\right)   &  \mapsto t;
\end{align*}
this is a genus one fibration since we have a cubic equation in $u,W.$

However, this fibration seems to have no section. Nevertheless, taking its
Jacobian fibration produces an elliptic fibration with section and the same
fiber type.

If we make a base change of this fibration: $\left(  t+s\right)  =m^{3}$ then
we obtain the following elliptic fibration with $U=um.$%
\[
U^{3}sm-\left(  s-m^{3}\right)  ksWU-W^{3}sm-\left(  s-m^{3}\right)
^{2}\left(  s^{2}-sm^{3}-1\right)  m=0.
\]
The transformation
\begin{align*}
W &  =\frac{-24y+12(s^{2}+1)(s-m^{3})x+\left(  s-m^{3}\right)  ^{2}Q}{18sm(4x-m^{2}(s^{2}+1)^{2})}\\
U &  =\frac{-24y-12(s^{2}+1)(s-m^{3})x+\left(  s-m^{3}\right)  ^{2}Q}{18sm(4x-m^{2}(s^{2}+1)^{2})}%
\end{align*}
where $Q=\left(
108m^{6}s^{3}+(s^{6}+105s^{2}-111s^{4}-1)m^{3}+s(s^{2}+1)^{3}\right)$
gives a Weierstrass equation of the Jacobian variety $J(Y_k^{(3)})$, the point $\pi_{3}$ of $x$ coordinate $\frac
{1}{4}m^{2}(s^{2}+1)^{2}$ being a $3$-torsion point. Taking again $m^{3}=\left(
t+s\right)  ,$  and $\pi_{3}=\left(  X=0,Y=0\right)  ,$ we recover a
Weierstrass equation for the $3$-isogenous fibration  $\#19$%

\[
Y^{2}-3tkYX-t^{2}\left(  27t^{2}-k\left(  k^{2}-27\right)  t+27\right)
Y=X^{3}%
\]
hence a fibration of $N_{k}.$
Recall that the transcendental lattice of $Y_k^{(3)}$ is $T(Y_k)[3]$ (Theorem \ref{th:sh} \cite{Sh1}. 
\end{proof}

\subsection{3-isogenies of $Y_2$}
\begin{theorem}
All the $3$-isogenous surfaces obtained from $3$-torsion sections of $Y_2$ are the $K3$-surface $Y_{10}$.

\end{theorem}
\begin{proof}
Recall the results, given in \cite{BL1}, about the $4$ elliptic fibrations of $Y_2$ 
with $3$-torsion. %

\[%
\begin{array}
[c]{ccc}%
&
\begin{array}
[c]{c}%
\text{Weierstrass Equation}\\
\text{Singular Fibers}%
\end{array}
& \text{Rank}\\
\hline
\#20\left(  7-w\right)   &
\begin{array}
[c]{c}%
E_w: Y^{2}-\left(  w^{2}+2\right)  YX-w^{2}Y=X^{3}\\
I_{12}\left(  \infty\right)  ,\quad I_{6}\left(  0\right)  ,\quad2I_{2}\left(
\pm1\right)  ,\quad2I_{1}%
\end{array}
& 0\\
\hline
\#19\left(  8-b\right)   &
\begin{array}
[c]{c}%
E_b: Y^{2}+2bYX+b^{2}\left(  b+1\right)  ^{2}Y=X^{3}\\
2IV^{\ast}\left(  \infty,0\right)  ,\quad I_{6}\left(  -1\right)  ,\quad2I_{1}%
\end{array}
& 1\\
\hline
20-j &
\begin{array}
[c]{c}%
E_j: Y^{2}-4\left(  j^{2}-1\right)  YX+4\left(  j+1\right)  ^{2}Y=X^{3}\\
I_{12}\left(  \infty\right)  ,\quad IV^{\ast}\left(  -1\right)  ,\quad
I_{2}\left(  -\frac{1}{2}\right)  ,\quad2I_{1}%
\end{array}
& 0\\
\hline
21-c &
\begin{array}
[c]{c}%
E_c:Y^{2}+\left(  c^{2}+5\right)  YX+Y=X^{3}\\
I_{18}\left(  \infty\right)  ,\quad6I_{1}%
\end{array}
& 1
\\\hline
\end{array}
\]

Using \ref{formulae} we compute the $3$-isogenous elliptic fibrations
named $H_{w},H_{b}$ $,H_{j}$ and $H_{c}$ and given in the next table. To simplify we denote $H_w$ (resp. $H_b$) the specialized elliptic fibration $H_{\#20}(2)$ (resp. $H_{\#19}(2)).$


\[
\begin{array}
[c]{cc}
\hline
&
\begin{array}
[c]{c}%
\text{Weierstrass Equation}\\
\text{Singular Fibers}%
\end{array}
\\
\hline
H_{w} &
\begin{array}
[c]{c}%
Y^{2}+3\left(  w^{2}+2\right)  YX+\left(  w^{2}+8\right)  \left(
w^{2}-1\right)  ^{2}Y=X^{3}\\
I_{4}\left(  \infty\right)  ,\quad2I_{6}\left(  \pm1\right)  ,\quad
I_{2}\left(  0\right)  ,\quad2I_{3}%
\end{array}
\\
\hline
H_{b} &
\begin{array}
[c]{c}%
Y^{2}-6bYX+b^{2}\left(  27b^{2}+46b+27\right)  Y=X^{3}\\
  2IV^{\ast}\left(  \infty,0\right)  ,\quad2I_{3},\quad I_{2}\left(  -1\right)%
\end{array}
\\
\hline
H_{j} &
\begin{array}
[c]{c}%
Y^{2}+12\left(  j^{2}-1\right)  YX+4\left(  4j^{2}-12j+11\right)  \left(
j+1\right)  ^{2}\left(  2j+1\right)  ^{2}Y=X^{3}\\
I_{4}\left(  \infty\right)  ,\quad IV^{\ast}\left(  -1\right)  ,\quad
I_{6}\left(  \frac{-1}{2}\right)  ,\quad2I_{3}%
\end{array}
\\
\hline
H_{c} & \begin{array}
[c]{c}%
Y^{2}-3\left(  c^{2}+5\right)  YX-\left(  c^{2}+2\right)  \left(
c^{2}+c +7\right)  \left(  c^{2}-c+7\right)  Y=X^{3}\\
I_{6}\left(  \infty\right)  ,\quad6I_{3}%
\end{array}
  \\
  \hline

  \end{array}
\]

We know from Theorem \ref{th:spe} that $H_w$ and $H_b$ are elliptic fibrations of $Y_{10}$;
we know also from the classification of  elliptic $K3$ surfaces of rank $0$ \cite{SZ} that $H_w$ and $H_j$, with rank $0$, are fibrations of $Y_{10}$.
Fibrations $H_j$ and $H_c$ are fibrations of $Y_{10}$ by applying theorem \ref{3isog}. From Shimada \cite{Shim2} we observe that on $H_c$ with fibers $I_6$ and $6I_3$ we get $(\mathbb Z /3\mathbb Z )^2$ torsion. We recover this last result using formula 2.7.2 giving the equation of $H_c$ as $u^3+v^3+w^3=(c^2+5)uvw$ thus showing that $H_c$ is an element of Hesse pencil. The nine inflexion points correspond to $3$-torsion sections.

 \begin{remark}
On the elliptic fibration $21-c$ there exists an involution $i_c$ $c \mapsto -c$, $X \mapsto X$, $Y_1=Y-(c^2+5)/2\mapsto -Y_1$ . We can show that 
there exists an elliptic fibration of parameter $X$ with a two torsion section, (fibration $24-\psi$ of $Y_2$ \cite{BL1} last table), such that $i_c$ corresponds to add the two-torsion section. The same remark is true for the fibration $8-b$. 
\end{remark}
\end{proof}

\subsection{$3$-isogenies of $Y_2$ and base-change of elliptic fibrations}

The $4$ fibrations $H_w$, $H_b$, $H_j$, $H_c$ of $Y_{10}$ obtained by $3$-isogenies have the same rank as the corresponding fibrations $E_w,E_b,E_j,E_c$ of $Y_2$.
On the contrary, in the previous theorem \ref{3isog} we study a base change $h=u^3$ of an elliptic fibration with 2 singular fibers of type $II^*$. This theorem applied with $A(t)=t^2+5$,$B(t)=1$ corresponding to the fibration $21-c$ of $Y_2$ \cite{BL1} and then fibration $13-h$ of rank $1$ defines an other elliptic fibration of $Y_{10}$ increasing the rank to $7.$ In this paragraph we give also two other examples of  base change of elliptic fibration of $Y_2$ leading to  elliptic fibrations of $Y_{10}$ of rank $4$.
We also construct generators of the Mordell-Weil
lattice of these fibrations.

With some elementary change of variables we use Weierstrass equations $E_h$, $E_f$, $E_g$ of respective fibrations $13-h,11-f,12-g$ of $Y_2$ given in \cite{BL1}.

\begin{theorem}\label{th:5.5}

    \begin{enumerate}
    \item The elliptic fibration $(E_h)$ of $Y_2$ with Weierstrass equation 
      \[(E_h) \qquad Y^2=X^3-5h^2X^2+h^5(h+1)^2,\] 
      and
      singular fibers $2II^*(0,\infty)$, $I_2(1)$, $2I_1$, gives by base change $h=u^3$ a rank $7$ elliptic fibration $(E_u)$
      of $Y_{10}$ with singular fibers 
$2I_{0}^{\ast}\left(  \infty,0\right)$, $3I_{2}\left(  u^{3}+1\right)$
, $6I_{1}.$  
    \item The elliptic fibration $(E_f)$ of $Y_2$ with Weierstrass equation 
      \[(E_f) \qquad Y^2=X^3-f(2f-3)X^2+3f^2(f-1)^2X+f^3(f-1)^4,\]
      and singular fibers $II ^*(\infty)$, $III^*(0)$, $I_4(1)$, $I_1(32/27)$, gives by base change $f=u'^3$ a rank $4$ elliptic fibration $(E_{u'})$ of $Y_{10}$ with singular fibers $I_0^*(\infty)$, $III(0)$, $3I_4(1,\zeta_3,\zeta_3^2)$, $3I_1$.

      The elliptic fibration $(E_g)$ of $Y_2$ with Weierstrass equation
      \[ (E_g) \qquad y^2=x^3+4g^2x^2+g^3(g+1)^2x,\]
      and singular fibers $2III^*(0,\infty)$, $I_4(-1)$, $I_2(1)$, gives by base change $g=n^3$ a rank $4$ elliptic fibration $(E_n)$ of $Y_{10}$ with singular fibers $2III(0,\infty)$, $3I_4(-1,n^2-n+1)$, $3I_2(1,n^2+n+1)$.

    \end{enumerate}      

  \end{theorem}

  \begin{proof}
    \begin{enumerate}

\item  For the first fibration $E_h$ with $h=u^3$ we get 
\begin{align*}
(E_{u}) &  :Y^{2}=X^{3}-5u^{2}X^{2}+u^{3}(u^{3}+1)^{2}\\
        &  2I_{0}^{\ast}\left(  \infty,0\right)  ,3I_{2}\left(  u^{3}+1\right) ,6I_{1}.
\end{align*}
This is a fibration of $Y_{10}$ from theorem \ref{3isog}.
Notice that the $j$-invariant of $(E_{u})$ is invariant by the two
transformations $u\mapsto\frac{1}{u}$ and $u\mapsto \zeta u$ ($\zeta^3=1$). These automorphisms
of the base $\mathbb{P}^{1}$ of the fibration $(E_{u})$ can be extended to the
sections as explained below.

Let $S_{3}=<\zeta_3,\tau;\zeta_3^{3}=1,\tau^{2}=1>$ be the non abelian group of
order $6$ and define an action of $S_{3}$ on the sections of $E_{u}$ by%
\begin{align*}
&  \left(  X(u),Y(u)\right)  \overset{\tau}{\mapsto}\left(  u^{4}X\left(
\frac{1}{u}\right)  ,u^{6}Y\left(  \frac{1}{u}\right)  \right)  \\
&  \left(  X\left(  u\right)  ,Y\left(  u\right)  \right)  \overset{\zeta_3
}{\mapsto}\left(  \zeta_3 X\left(  \zeta_3 u\right)  ,Y\left(  \zeta_3 u\right)  \right).
\end{align*}

To obtain generators of $(E_{u})$, following Shioda \cite{Sh1}, we use the rational
elliptic surface $X^{\left(  3\right)  +}$ with $\sigma=u+\frac{1}{u}$ and a
Weierstrass equation%
\begin{align*}
(E_{\sigma})  &  :y^{2}=x^{3}-5x^{2}+\left(  \sigma-1\right)  ^{2}\left(
\sigma+2\right) \\
&  I_{0}^{\ast}\left(  \infty\right)  ,I_{2}\left(  -1\right)  ,4I_{1}%
\end{align*}
of rank $3.$
The Mordell-Weil lattice of a rational elliptic surface is generated by
sections of the form $(a+b\sigma+c\sigma^{2},d+e\sigma+f\sigma^{2}+g\sigma
^{3}).$ Moreover since we have a singular fiber of type $I_{0}^{\ast}$ at
$\infty$ the coefficients $c$ and $f,g$ are $0$ \cite{KK}. So after an easy
computation we find the $3$ sections \ (with $\zeta_3^{3}=1,i^{2}=-1).$%
\begin{align*}
q_{1}  & =\left(  -\left(  \sigma-1\right)  ,i\sqrt{2}\left(  \sigma-1\right)
\right)  \quad\\
q_{2}  & =\left(  -\zeta_3 \left(  \sigma-1\right)  ,\left(  3+\zeta_3 \right)  \left(
\sigma-1\right)  \right)  \quad q_{3}=\left(  -\zeta_3^{2}\left(  \sigma-1\right)
,\left(  3+\zeta_3^{2}\right)  \left(  \sigma-1\right)  \right).
\end{align*}
These sections give the sections $\pi_{1},\pi_{2},\pi_{3}$ on $(E_{u})$ which
are fixed by $\tau$.
\[%
\begin{array}
[c]{c}%
\pi_{1}=\left(  -u\left(  u^{2}-u+1\right)  ,i\sqrt{2}u^{2}\left(
u^{2}-u+1\right)  \right)  \\
\pi_{2}=\left(  -\zeta_3 u\left(  u^{2}-u+1\right)  ,\left(  3+\zeta_3 \right)  u^{2}\left(
u^{2}-u+1\right)  \right)  \\
\pi_{3}=\left(  -\zeta_3^{2}u\left(  u^{2}-u+1\right)  ,\left(  3+\zeta_3^{2}\right)
u^{2}\left(  u^{2}-u+1\right)  \right).
\end{array}
\]  

We notice $\rho_{i}=\gamma\left(  \pi_{i}\right)  $ and $\mu_{i}=\gamma
^{2}\left(  \pi_{i}\right)  $ for $1\leq i\leq3$ which give $9$ rational
sections with some relations.

Moreover we have another section from the fibration $E_{h}$ of rank
$1.$

The point of $x$ coordinate $\frac{1}{16}\left(  h^{2}+\frac{1}{h^{2}}\right)
-h-\frac{1}{h}+\frac{29}{24}$ is defined on $\mathbb{Q}\left(  h\right)  .$
Passing to $(E_{u})$ we obtain $\omega=$
\[
\left(  \frac{1}{16}\frac{(  1-16u^{3}+46u^{6}-16u^{9}%
+u^{12})  }{u^{4}},-\frac{1}{64}\frac{(  u^{6}-1)  (
1-24u^{3}+126u^{6}-24u^{9}+1)  }{u^{6}}\right ).
\]

We hope to get a generator system with $\pi_{i},\rho_{i}$ and $\omega$ so we
have to compute the height-matrix. The absolute value of its determinant  is
$\frac{81}{16}.$ Since the discriminant of the surface is $72,$ we obtain a
subgroup of index $a$ with $\frac{81}{16}\times\frac{1}{a^{2}}\times2^{3}%
4^{2}=72$ so $a=3$.

After some specialization of $u$ $\in\mathbb{Z}$ (for example if $u=11$, $(E_{u})$
has rank $3$ on $\mathbb{Q})$ we find other sections with $x$ coordinate
of the shape \ $\ \left(  au+b\right)  \left(  u^{2}-u+1\right)  $
\begin{align*}
\mu &  =\left(  -\left(  u-1\right)  \left(  u^{2}-u+1\right)  ,-\left(
u^{2}-u+1\right)  \left(  u^{2}+2u-1\right)  \right)  \\
\mu_{1} &  =\left(  -\left(  u-9\right)  \left(  u^{2}-u+1\right)  ,\left(
u^{2}-u+1\right)  \left(  5u^{2}-18u+27\right)  \right)  \\
\mu_{2} &  =\left(  -\left(  u+\frac{1}{3}\right)  \left(  u^{2}-u+1\right)
,\frac{i\sqrt{3}}{9}\left(  u^{2}-u+1\right)  \left(  9u^{2}+4u+1\right)
\right).
          \end{align*}
We deduce the relations
\begin{align*}
3\mu &  =\omega+\pi_{2}-\gamma\left(  \pi_{3}\right)  +\pi_{3}-\gamma
^{2}\left(  \pi_{2}\right)  \\
&  =\omega+2\pi_{2}+\gamma\left(  \pi_{2}\right)  +\pi_{3}-\gamma\left(
\pi_{3}\right)
\end{align*}
so, the Mordell-Weil lattice is generated by $\pi_{j},\rho_{j}=\gamma\left(
\pi_{j}\right)  $ for $1\leq j\leq3$ and $\mu$ with the following Gram matrix of determinant $9/16$ (\cite{Sh3} Corollary 1.7).%

\[
\left(
\begin{array}
[c]{ccccccc}%
1 & -\frac{1}{2} & 0 & 0 & 0 & 0 & 0\\
-\frac{1}{2} & 1 & 0 & 0 & 0 & 0 & 0\\
0 & 0 & 1 & -\frac{1}{2} & 0 & 0 & \frac{1}{2}\\
0 & 0 & -\frac{1}{2} & 1 & 0 & 0 & 0\\
0 & 0 & 0 & 0 & 1 & -\frac{1}{2} & \frac{1}{2}\\
0 & 0 & 0 & 0 & -\frac{1}{2} & 1 & \frac{-1}{2}\\
0 & 0 & \frac{1}{2} & 0 & \frac{1}{2} & \frac{-1}{2} & 2
\end{array}
\right)  .
\]

\item
  
We consider now the second and third fibration, respectively $E_f$ and $E_g$.

First we explain why such base changes lead to elliptic fibrations of $Y_{10}.$

Recall (see \cite{BL1}) that, from the equation
\[ Y_{2}:X+\frac{1}{X}+Y+\frac{1}{Y}+Z+\frac{1}{Z}=2
\]
the  fibration with two singular fibers $II^{\ast}$ can be obtained with the
parameter $h=\frac{\left(  Y+Z\right)  YX^{2}}{Z^{3}\left(  X+Z\right)  }.$ We
can also obtain with the parameter $g=\frac{Y\left(  Y+Z\right)  }{X\left(
X+Z\right)  }$ a fibration with singular fibers $2III^{\ast},I_{4}(0,\infty)$ and
$I_{2}$ (fibration $\left(  12-g\right)  $ of $Y_2$ (\cite{BL1} last table)). It is also possible to get with the parameter $f=\frac{Y\left(
X+Z\right)  ^{2}\left(  Z+Y\right)  }{Z^{3}X}$ singular fibers $II^{\ast
},III^{\ast},I_{4}$ and $I_{1}$(fibration $\left(  11-f\right)  $ of $Y_{2}$(\cite{BL1} last table)). We observe that $g=\frac{Z^{3}}{X^{3}}h,$ and
$f=\frac{\left(  X+Z\right)  ^{3}}{Z^{3}}h.$ So since the fibration of $Y_{10}$
obtained by a change of basis $h=u^{3}$ of the corresponding fibration allow
also to construct with a change of bases $g=r^{3}$ and $f=u^{\prime3}$ two other
fibrations of the same surface $Y_{10}$. 

We give some details of how this $2$ fibrations are linked to fibrations of $Y_2$ with $3$-torsion sections and also their Mordell-Weil generators.

\begin{enumerate}
  \item For the parameter $f$ we start from Weierstrass equation $E_w$ in the previous table and consider an other fibration with the parameter $f=Y$. Then we
obtain the fibration $E_f$  with singular fibers $II^{\ast
}\left(  \infty\right)  ,III^{\ast}\left(  0\right)  ,$ $I_{4}\left(  1\right)
,I_{1}\left(  \frac{32}{27}\right)  .$ A Weierstrass equation
\[
E_{f}:Y_{1}^{2}=X_{1}^{3}-f(2f-3)X_{1}^{2}+3f^{2}\left(  f-1\right)  ^{2}%
X_{1}+f^{3}\left(  f-1\right)  ^{4}%
\]
is obtained from $E_w$ with the following transformations%
\begin{align*}
f  &  =Y,\quad X_{1}=-\frac{Y\left(  Y-1\right)  ^{2}}{X+1},\quad Y_{1}%
=w\frac{Y^{2}\left(  Y-1\right)  ^{2}}{X+1}\\
w  &  =\frac{-Y_{1}}{X_{1}f},\quad X=-\frac{X_{1}+f\left(  f-1\right)  ^{2}%
}{X_{1}},\quad Y=f.
\end{align*}
The base change $f(=Y)=u^{\prime3} $  gives an elliptic
fibration of $Y_{10}$ with singular fibers $I_{0}^{\ast}\left(  \infty\right)
$,  $III\left(  0\right)  $,  $3I_{4}\left(  1,\zeta,\zeta^{2}\right)  $,  $3I_{1}$ , rank $4$,
a Weierstrass equation and sections
\begin{align*}
{y^{\prime}}^2  &  =x^{\prime3}+u^{\prime2}x^{\prime2}+2u^{\prime}\left(
u^{\prime3}-1\right)  x^{\prime}+u^{\prime3}\left(  u^{\prime3}-1\right)  ^{2}\\
P  &  =\left(  x_{P}\left(  u^{\prime}\right)
,y_{P}\left(  u^{\prime}\right)  \right)  =\left(  -\left(  u^{\prime
3}-1\right)  ,\left(  u^{\prime}-1\right)  ^{2}\left(  u^{\prime2}+u^{\prime
}+1\right)  \right) \\
Q  &  =\left(  x_{Q}\left(  u^{\prime}\right)  ,y_{Q}\left(  u^{\prime
}\right)  \right)  =\left(  -\left(  u^{\prime}+2\right)  \left(  u^{\prime
2}+u^{\prime}+1\right)  ,2i\sqrt{2}\left(  u^{\prime2}+u^{\prime}+1\right)
^{2}\right) . \label{rk41}
\end{align*}
Using also the points $P^{\prime}$ with $x_{P^{\prime}}=\zeta_3 x_{P}\left(
\zeta_3 u^{\prime}\right)  $ and $Q^{\prime}$ with $x_{Q^{\prime}}=\zeta_3 x_{Q}\left(
\zeta_3 u^{\prime}\right)  $ (with $\zeta_3^{3}=1$), we can
compute the height matrix and show that the Mordell-Weil  lattice is generated
by $P,P^{\prime},Q,Q^{\prime}$  and is equal to $A_2\left(  \frac{1}%
{4}\right)  \oplus A_2\left( \frac{1}{2}\right).$

\item 
The second example is obtained from $E_w$ with the parameter 
$g=\frac{Y}{w^{2}}$. We obtain  the fibration
$E_g$ with singular fibers $2III^*(0,\infty),I_4(-1),I_2(1)$
\begin{align*}
E_{g}  & :y^{2}=x^{3}+4x^{2}g^{2}+g^{3}(g+1)^{2}x
\end{align*}

with the following transformation
\begin{align*}
X  & =\frac{x^{2}\left(  x+2g^{2}\right)  \left(  g-1\right)  }{y^{2}}%
,Y=\frac{1}{g}\frac{x^{2}\left(  x+2g^{2}\right)  ^{2}}{y^{2}},w=\frac
{x\left(  x+2g^{2}\right)  }{gy}\\
x  & =\frac{Y^{2}\left(  Y-2X-w^{2}\right)  }{Xw^{4}},y=\frac{Y^{3}\left(
Y-w^{2}\right)  \left(  Y-2X-w^{2}\right)  }{X^{2}w^{7}},g=\frac{Y}{w^{2}}.%
\end{align*}
So in the Weierstrass equation $E_{g}$, if we replace the parameter $g(=Y/w^2)$ by
$n^{3}$ we obtain the following fibration of $Y_{10}$%
\label{rk42}\begin{align}
y^{2}  & =x^{3}+4n^{2}x^{2}+n\left(  n+1\right)^{2}\left(  n^{2}-n+1\right)^{2}x
\end{align}
with singular fibers $2III\left(  0,\infty\right)  ,3I_{4}\left(
-1,n^{2}-n+1\right)  ,3I_{2}\left(  1,n^{2}+n+1\right)  $ and rank $4.$
Notice the two infinite sections with $x$ coordinates $(n+1)^2(n^2-n+1)$ and $-\frac{1}{3}(n-1)^2(n^2-n+1).$
The previous method for the Mordell-Weil generators gives no results but since we have a "self 2-isogeny"
 we can use the complex multiplication to obtain $4$ generators. 

\end{enumerate}
\end{enumerate}
\end{proof}

\subsection{A fibration for Theorem \ref{2iso}}
From the previous fibration $E_u$ we construct by a $2$-neigbour method a fibration with a $2$-torsion section used in Theorem \ref{2iso}.

We start from the Weierstrass equation
\[
Y^{2}=X^{3}-5u^{2}X^{2}+u^{3}\left(  u^{3}+1\right)  ^{2}%
\]
and obtain another elliptic fibration with the parameter $m=\frac{X}{u\left(
u^{2}-u+1\right)  }$, which gives the Weierstrass equation
\begin{align}
E_{m}  &  :y^{2}=x^{3}-\left(  m^{3}+5m^{2}-2\right)  x^{2}+\left(
m^{3}+1\right)  ^{2}x \label{rk4}
\end{align}
with singular fibers 
$  I_{0}^{\ast}\left(  \infty\right) $, $3I_{4}\left(  m^{3}+1\right)
$,  $ I_{2}\left(  0\right)  $,  $4I_{1}\left(  1,-\frac{5}{3}%
,m^{2}-4m-4\right)$ 
 and rank $4.$ 
\begin{remark}
From this fibration with the parameter $q=\frac{y}{xm}$ we recover the
fibration $H_{c.}$
\end{remark}

\section{Some $3$-isogenies from elliptic fibrations of $Y_{10}$}
We found, along our study of elliptic fibrations of the $K3$ surface $Y_{10}$, some elliptic fibrations with $3$-torsion section. A rank $0$ (extremal) elliptic fibration numbered $80$ in Shimada-Zhang \cite{SZ} together with a rank $2$ elliptic fibration $(11)$ are given in Bertin and Lecacheux's paper \cite{BL3}. In the present paper, apart from the two specialized $3$-torsion elliptic fibrations of $Y_{10}$, namely $H_{\#19}(10)$ and  $H_{\#20}(10)$, we get the elliptic fibration $H_c$ with $\mathbb Z/3\mathbb Z \times \mathbb Z/3 \mathbb Z$-torsion. We are going to compute the transcendental lattices of $K3$ surfaces obtained by $3$-isogeny.
\begin{theorem}
  \begin{enumerate}
  \item Elliptic fibrations numbered $80$ (See Shimada's notation \cite{SZ}) and $(11)$ \cite{BL3} induce by $3$-isogeny an elliptic fibration of a $K3$ surface with transcendental lattice $[4 \quad 0 \quad 18]$.

  \item There is an elliptic fibration of $Y_{10}$ with a $3$-torsion section whose $3$-isogenous surface is a $K3$ surface with transcendental lattice $[2 \quad 0 \quad 36]$.

    \end{enumerate}

  \end{theorem}

  \begin{proof}

   Recall first the corresponding Weierstrass equations.

\[
\begin{array}
[c]{ccc}%
\hline
&
\begin{array}
[c]{c}%
\text{Weierstrass Equation}\\
\text{Singular Fibers}%
\end{array}
& \text{Rank}\\
\hline
\begin{array}
[c]{c}%
80\\
\text{Shimada notation}%
\end{array}
&
\begin{array}
[c]{c}%
Y^{2}+\left(  9t^{2}+6t-9\right)  YX+9t^{4}\left(  3t^{2}+6t-5\right)
Y=X^{3}\\
I_{12}\left(  0\right)  ,3I_{3}\left(  \infty,3t^{2}+6t-5\right)
,I_{2}\left(  \frac{3}{5}\right)  ,I_{1}\left(  -\frac{3}{4}\right)
\end{array}
& 0\\
\hline
E_{11} &
\begin{array}
[c]{c}%
Y^{2}+\left(  t^{2}-4\right)  YX+t^{2}\left(  2t^{2}-3\right)  Y=X^{3}\\
2I_{6}\left(  0,\infty\right)  ,2I_{3}\left(  2t^{2}-3\right)  ,2I_{2}\left(
\pm1\right)  ,2I_{1}\left(  \pm8\right)
\end{array}
& 2\\
\hline
\#19 &
\begin{array}
[c]{c}%
Y^{2}+10tYX+t^{2}\left(  t^{2}+10t+1\right)  Y=X^{3}\\
2IV^{\ast}\left(  0,\infty\right)  ,2I_{3}\left(  t^{3}+10t+1\right)
,2I_{1}\left(  \left(  t-27\right)  \left(  27t-1\right)  \right)
\end{array}
& 2\\
\hline
\#20 &
\begin{array}
[c]{c}%
Y^{2}-(t^{2}-10t+3)YX-\left(  t^{2}-10t+1\right)  Y=X^{3}\\
I_{12}\left(  \infty\right)  ,2I_{3}\left(  t^{2}-10t+1\right)  ,2I_{2}\left(
0,10\right)  ,2I_{1}\left(  \left(  t-1\right)  \left(  t-9\right)  \right)
\end{array}
& 1\\
\hline
H_{c} &
\begin{array}
[c]{c}%
Y^{2}-\left(  t^{2}+11\right)  YX-\left(  t^{2}+t+7\right)  \left(
t^{2}-t+7\right)  Y=X^{3}\\
I_{6}\left(  \infty\right)  ,6I_{3}\left(  t^{2}+2,t^{2}\pm t+7\right)
\end{array}
& 1\\
\hline
\end{array}
\]
Notice that from Theorem \ref{th:spe} the $3$-isogenous $K3$-surfaces of $H_{\#19}$ and$H_{\#20}$ have the lattice $[4 \quad 0 \quad 18]$ as transcendental lattice.

\begin{enumerate}
   
 \item
 \textbf{Fibration $80$ of rank $0$}

A Weierstrass equation for the $3$-isogenous fibration is 
\begin{align*}
 Y^{2}+\left(  -27t^{2}-18t+27\right)  YX+27\left(  4t+3\right)  \left(
5t-3\right)  ^{2}Y=X^{3}
\end{align*}
with singular fibers $ I_{9}\left(  \infty\right) ,I_{6}\left(  \frac{3}{5}\right),I_{4},\left(  0\right) ,I_{3}\left(  -\frac{3}{4}\right),2I_1.$

From singular fibers, torsion and rank we see in Shimada-Zhang table \cite{SZ} that it is the number $48$ case with transcendental lattice $[4 \quad 0 \quad 18]$.

\textbf{Fibration (11)}

  In section $6$ of Bertin and Lecacheux's paper \cite{BL3}  is exhibited a rank $2$ elliptic fibration named $(11)$ with Weierstrass equation
  \[y^2+(t^2-4)xy-t^2(t^2-63)y=(x-9t^2)(x^2+108t^2).\]

After a translation to put the $3$-torsion section in $\left(  0,0\right)  $
we obtain the  Weierstrass equation $E_{11}$ given in the table and also
$p_{1}$ and $p_{2}$ generators of the Mordell-Weil lattice
\begin{equation}%
\begin{array}
[c]{c}%
p_{1}=\left(  6t^{2},27t^{2}\right)  ,p_{2}=\left(  6i\sqrt{3}t-3t^{2}%
,27t^{2}\right)  \\
2I_{6}\left(  \infty,0\right)  ,\quad2I_{3}\left(  2t^{2}-3\right)
,\quad2I_{2}\left(  \pm1\right)  ,\quad2I_{1}\left(  \pm8\right).
\end{array}
\end{equation}
We deduce the Weierstrass equation of its $3$-isogenous elliptic fibration and its singular fibers
\begin{align*}
H_{11}: &  Y^{2}-3\left(  t^{2}-4\right)  YX-\left(  t^{2}-1\right)
^{2}\left(  t^{2}-64\right)  Y=X^{3}\\
&  2I_{6}\left(  \pm1\right)  ,\quad2I_{3}\left(  \pm8\right)  ,\quad
2I_{2}\left(  \infty,0\right)  ,\quad2I_{1}\left(  2t^{2}-3\right)  \text{ of
rank }2.
\end{align*}
Moreover, computing the height matrix of the two sections
\begin{align*}
\pi_{1}  & =\left(  -\frac{1}{4}\left(  t^{2}-1\right)  \left(  t^{2}%
-64\right)  ,\frac{1}{8}\left(  t-8\right)  \left(  t-1\right)  \left(
t+1\right)  ^{2}\left(  t+8\right)  ^{2}\right)  \quad\\
\omega & =\left(  -7\left(  t^{2}-1\right)  ^{2},49\alpha\left(
t^{2}-1\right)  ^{3}\right)  \quad\text{where }49\alpha^{2}+20\alpha+7=0,
\end{align*}
we find that the discriminant of the $K3$ surface is $72.$

Let us compute the transcendental lattice of $H_{11}$.

For each reducible fiber at $t=i$ we denote $\left(  X_{i},Y_{i}\right)  $ the
singular point of $H_{11}$
\[%
\begin{array}
[c]{ccc}%
t=\pm1 \quad t=\pm8 & t=0 & t=\infty\\
(  X_{\pm1}=0,Y_{\pm1}=0) & 
(  X_{0}=-16,Y_{0}=64)  & (  x_{\infty}=-1,y_{\infty
}=-1)\\
(  X_{\pm8}=0,Y\pm_{8}=0)& &
\end{array}
\]
where if $t=\infty$ we substitute $t =\frac{1}{T}$,
$x=T^{4}X,y=T^{6}Y.$ We denote also $\theta_{i,j}$ the $j$-th component
of the reducible fiber at $t=i.$ A section $M=\left(  X_{M},Y_{M}\right)  $
intersects the component $\theta_{i,0}$ if and only if $\left(  X_{M},Y_{M}\right)
\not \equiv\left(  X_{i},Y_{i}\right)  \operatorname{mod}\left(  t-i\right)
.$ Using the additivity on the component, we deduce that $\omega$ does not
intersect $\theta_{i,0},2\omega$ intersects $\theta_{i,0}$ and so $\omega$
intersects $\theta_{i,3}$ for $i=\pm1.$ Also $\omega$ intersects $\theta_{i,0}$
for $i=\pm8$ $,i=0$ and $\infty.$

For $\pi_{1}$ we compute $k\pi_{1}$ with $2\leq k\leq6.$ For $i=\pm1$, only
$6\pi_{1}$ intersects $\theta_{i,0}$ so $\pi_{1}$ intersects $\theta_{i,1}$
(this choice 1, not 5, gives the numbering of components). For $i=\pm8$, only
$3\pi_{1}$ intersects $\theta_{i,0},$ so $\pi_{1
}$ intersects $\theta_{i,1}.$
Modulo $t$,\ we get $\pi_{1}=\left(  -16,64\right)  ,$ so $\pi_{1}$ intersects
$\theta_{0,1},$ and $\pi_{1}$ intersects $\theta_{\infty,0}$.

As for the $3$-torsion section $s_{3}=\left(  0,0\right)  $,  $s_{3}$
intersects $\theta_{i,2}$ or $\theta_{i,4}$ if $i=\pm1$. Computing $2\pi_{1}-s_{3}$, it follows that $s_{3}$ intersects $\theta_{1,2}$ and $\theta_{-1,4}$.

For $i=\pm8$, we compute $\pi_{1}-s_{3},$ for $t=8$ and show that $s_{3}$
intersects $\theta_{8,2} $ and
$\theta_{-8,1}.$ For $t=0$ and $t=\infty\ s_{3}$ intersects \ the $0$ component.

From the relation 
$\theta_{i,j}$ and find  that $3s_{3}\approx -2\theta_{1,1}-4%
\theta_{1,2}-3\theta_{1,3}-2\theta_{1,4}-\theta
_{1,5}$, we can choose the following base of the N\'{e}ron-Severi lattice ordered as 
$s_{0},F,\theta_{1,j},$ with $1\leq j\leq4,$ $s_{3},\theta_{-1,k}$ with $1\leq
k\leq5,\theta_{8,k},k=1,2,$ $\theta_{-8,k},k=1,2$ and $\theta_{0,1},$
$\theta_{\infty,1},\omega,\pi_{1}.$

Finally only the two sections $\omega$ and $\pi_{1}$
intersect. Hence it follows the Gram matrix $Ne$ of the N\'{e}ron-Severi
lattice,
\[Ne=
\left(
\begin{smallmatrix}
-2 & 1 & 0 & 0 & 0 & 0 & 0 & 0 & 0 & 0 & 0 & 0 & 0 & 0 & 0 & 0 & 0 & 0 & 0 &
0\\
1 & 0 & 0 & 0 & 0 & 0 & 1 & 0 & 0 & 0 & 0 & 0 & 0 & 0 & 0 & 0 & 0 & 0 & 1 &
1\\
0 & 0 & -2 & 1 & 0 & 0 & 0 & 0 & 0 & 0 & 0 & 0 & 0 & 0 & 0 & 0 & 0 & 0 & 0 &
1\\
0 & 0 & 1 & -2 & 1 & 0 & 1 & 0 & 0 & 0 & 0 & 0 & 0 & 0 & 0 & 0 & 0 & 0 & 0 &
0\\
0 & 0 & 0 & 1 & -2 & 1 & 0 & 0 & 0 & 0 & 0 & 0 & 0 & 0 & 0 & 0 & 0 & 0 & 1 &
0\\
0 & 0 & 0 & 0 & 1 & -2 & 0 & 0 & 0 & 0 & 0 & 0 & 0 & 0 & 0 & 0 & 0 & 0 & 0 &
0\\
0 & 1 & 0 & 1 & 0 & 0 & -2 & 0 & 0 & 0 & 1 & 0 & 0 & 1 & 1 & 0 & 0 & 0 & 0 &
0\\
0 & 0 & 0 & 0 & 0 & 0 & 0 & -2 & 1 & 0 & 0 & 0 & 0 & 0 & 0 & 0 & 0 & 0 & 0 &
1\\
0 & 0 & 0 & 0 & 0 & 0 & 0 & 1 & -2 & 1 & 0 & 0 & 0 & 0 & 0 & 0 & 0 & 0 & 0 &
0\\
0 & 0 & 0 & 0 & 0 & 0 & 0 & 0 & 1 & -2 & 1 & 0 & 0 & 0 & 0 & 0 & 0 & 0 & 1 &
0\\
0 & 0 & 0 & 0 & 0 & 0 & 1 & 0 & 0 & 1 & -2 & 1 & 0 & 0 & 0 & 0 & 0 & 0 & 0 &
0\\
0 & 0 & 0 & 0 & 0 & 0 & 0 & 0 & 0 & 0 & 1 & -2 & 0 & 0 & 0 & 0 & 0 & 0 & 0 &
0\\
0 & 0 & 0 & 0 & 0 & 0 & 0 & 0 & 0 & 0 & 0 & 0 & -2 & 1 & 0 & 0 & 0 & 0 & 0 &
1\\
0 & 0 & 0 & 0 & 0 & 0 & 1 & 0 & 0 & 0 & 0 & 0 & 1 & -2 & 0 & 0 & 0 & 0 & 0 &
0\\
0 & 0 & 0 & 0 & 0 & 0 & 1 & 0 & 0 & 0 & 0 & 0 & 0 & 0 & -2 & 1 & 0 & 0 & 0 &
1\\
0 & 0 & 0 & 0 & 0 & 0 & 0 & 0 & 0 & 0 & 0 & 0 & 0 & 0 & 1 & -2 & 0 & 0 & 0 &
0\\
0 & 0 & 0 & 0 & 0 & 0 & 0 & 0 & 0 & 0 & 0 & 0 & 0 & 0 & 0 & 0 & -2 & 0 & 0 &
1\\
0 & 0 & 0 & 0 & 0 & 0 & 0 & 0 & 0 & 0 & 0 & 0 & 0 & 0 & 0 & 0 & 0 & -2 & 0 &
0\\
0 & 1 & 0 & 0 & 1 & 0 & 0 & 0 & 0 & 1 & 0 & 0 & 0 & 0 & 0 & 0 & 0 & 0 & -2 &
1\\
0 & 1 & 1 & 0 & 0 & 0 & 0 & 1 & 0 & 0 & 0 & 0 & 1 & 0 & 1 & 0 & 1 & 0 & 1 & -2
\end{smallmatrix}
\right ).
\]

According to Shimada's lemma \ref{lem:gram},  $G_{Ne}\equiv \mathbb{Z/}2\mathbb{Z} \oplus\mathbb{Z/}36\mathbb{Z}$ is generated by the vectors $L_1$ and $L_2$ satisfying $q_{L_{1}}=-\frac{1}{2},$ $q_{L_{2}}=\frac{37}{36},$ and $b(L_1,L_2)=\frac{1}{2}.$

Since the following generators of the discriminant group of the lattice with Gram matrix 
$M_{18}=\left(
\begin{smallmatrix}
4 & 0\\
0 & 18
\end{smallmatrix}
\right)$ namely $f_1=(0,\frac{1}{2})$, $f_2=(\frac{1}{4},\frac{7}{18})$ verify $q_{f_1}=\frac{1}{2}$, $q_{f_2}=-\frac{37}{36}$ and $b(f_1,f_2)=-\frac{1}{2}$, we deduce that the Gram matrix of the transcendental lattice of the surface is $[4 \quad 0 \quad 18]$.

\item
\textbf{Fibration $H_c$}

We have shown along the proof of Theorem 5.4 that $H_c$ has a $(\mathbb Z/3\mathbb Z)^2$- torsion group 
and exhibited a 3-isogeny between some elliptic fibrations of $Y_{10}$ and $Y_2$. 
In the previous Weierstrass equation $H_{c}$ the point
$\sigma_{3}$ of $X$ coordinate $-\left(  c^{2}+c+7\right)  \left(
  c^{2}-c+7\right)  $ defines a $3$-torsion section. Notice that the group $\langle \sigma_3\rangle$ is stable by the Galois group $Gal(\bar{\mathbb Q}/\mathbb Q)$.  

After a translation to put this point in $(0,0)$ 
and scaling, we obtain a Weierstrass equation 
\[
Y'^2-(c^2+11)Y'X'-(c^2+c+7)(c^2-c+7)Y'=X'^3
.\]

The $3$-isogenous curve of
kernel $<\sigma_{3}>$ has a Weierstrass equation  %
\begin{align*}
y^{2}+3\left(  c^{2}+11\right)  xy+\left(  c^{2}+2\right)  ^{3}y &  =x^{3},
\end{align*}
with singular fibers $I_{2}\left(  \infty\right)  $, $2I_{9}\left(
  c^{2}+2\right)  $, $4I_{1}\left(  c^{4}+13c^{2}+49\right)  .$

The section
$P_{c}=\left(  -\frac{1}{4}\left(  c^{4}+c^{2}+1\right)  ,-\frac{1}{8}\left(
c-\zeta_3\right)  ^{3}\left(  c+\zeta_3^{2}\right)  ^{3}\right)  $ of infinite order, generates the Mordell-Weil lattice.

We consider the components of the reducible fibers in the following order
$\theta_{i\sqrt{2},j}$ $j\leq1\leq8,$ $\theta_{-i\sqrt{2},k}$ $1\leq k\leq8$,
$\theta_{\infty,1}$.

The $3$-torsion section $s_{3}=\left(  0,0\right)  $ and the previous
components are linked by the relation
\[
3s_{3}\approx -\theta_{i\sqrt{2},8}+\sum a_{i,j}\theta_{i,j}.%
\]
So we can replace, in the previous ordered sequence of components, the element
$\theta_{i\sqrt{2},8}$ by $s_{3}.$ We notice that $\left(  s_{3}.P_{c}\right)
=2,\left(  P_{c}.\theta_{\pm i\sqrt{2},0}\right)  =1$ and $\left(
P_{c}.\theta_{\infty,0}\right)  =1,$ so  the Gram matrix of the N\'{e}ron-Severi lattice is
\[
\left (\begin {smallmatrix} -2&1&0&0&0&0&0&0&0&0&0&0&0
&0&0&0&0&0&0&0\\1&0&0&0&0&0&0&0&0&1&0&0&0&0&0&0&0&0&0
&1\\0&0&-2&1&0&0&0&0&0&0&0&0&0&0&0&0&0&0&0&0
\\0&0&1&-2&1&0&0&0&0&0&0&0&0&0&0&0&0&0&0&0
\\0&0&0&1&-2&1&0&0&0&1&0&0&0&0&0&0&0&0&0&0
\\0&0&0&0&1&-2&1&0&0&0&0&0&0&0&0&0&0&0&0&0
\\0&0&0&0&0&1&-2&1&0&0&0&0&0&0&0&0&0&0&0&0
\\0&0&0&0&0&0&1&-2&1&0&0&0&0&0&0&0&0&0&0&0
\\0&0&0&0&0&0&0&1&-2&0&0&0&0&0&0&0&0&0&0&0
\\0&1&0&0&1&0&0&0&0&-2&0&0&1&0&0&0&0&0&0&2
\\0&0&0&0&0&0&0&0&0&0&-2&1&0&0&0&0&0&0&0&0
\\0&0&0&0&0&0&0&0&0&0&1&-2&1&0&0&0&0&0&0&0
\\0&0&0&0&0&0&0&0&0&1&0&1&-2&1&0&0&0&0&0&0
\\0&0&0&0&0&0&0&0&0&0&0&0&1&-2&1&0&0&0&0&0
\\0&0&0&0&0&0&0&0&0&0&0&0&0&1&-2&1&0&0&0&0
\\0&0&0&0&0&0&0&0&0&0&0&0&0&0&1&-2&1&0&0&0
\\0&0&0&0&0&0&0&0&0&0&0&0&0&0&0&1&-2&1&0&0
\\0&0&0&0&0&0&0&0&0&0&0&0&0&0&0&0&1&-2&0&0
\\0&0&0&0&0&0&0&0&0&0&0&0&0&0&0&0&0&0&-2&0
\\0&1&0&0&0&0&0&0&0&2&0&0&0&0&0&0&0&0&0&-2
\end {smallmatrix}\right).
\]
Its determinant is $-72.$
According to Shimada's lemma \ref{lem:gram}, $G_{NS}=\mathbb{Z}/2\mathbb{Z}%
\oplus\mathbb{Z}/36$ is generated by the vectors $L_{1}$ and $L_2$ satisfying 
$q_{L_1}=\frac{-1}{2},q_{L_2}=\frac{5}{36}$ and $b(L_1,L_2)=\frac{1}{2}.$

Moreover the following generators of the discriminant group of the lattice with Gram matrix
$M_{36}=\left(\begin{smallmatrix}
2&0\\
0&36
\end{smallmatrix}\right)$, namely $f_1=(\frac{1}{2},0)$, $f_2=(\frac{-1}{2},\frac{-7}{36})$ verify $q_{f_1}=\frac{1}{2}$, $q_{f_2}=-\frac{5}{36}$ and 
$b(f_1,f_2)=\frac{-1}{2}.$
So the Gram matrix of the transcendental lattice is $[2 \quad 0 \quad 36]$.

\end{enumerate}

\end{proof}

\section{Isogenies as isometries of the rational transcendental lattice}

Denoting the rational transcendental lattice $T(X)_{\mathbb Q}:=T(X) \otimes \mathbb Q$, we recall that $T(X)_{\mathbb Q}$ and $T(Y)_{\mathbb Q}$ are isometric if they define congruent lattices, that is if there exists $M\in Gl_n(\mathbb Q)$ satisfying $T(X)_{\mathbb Q}=^tMT(Y)_{\mathbb Q}M$.

Boissi\`ere, Sarti and Veniani proved the following theorem \cite{BSV}

\begin{theorem} \cite{BSV}
Let $\gamma: X \rightarrow Y$ be a $p$-isogeny between complex projective $K3$ surfaces $X$ and $Y$. Then $\text{rk}( T(Y)_{\mathbb {Q}})=\text{rk}( T(X)_{\mathbb {Q}})=:r$ and
\begin{enumerate}
\item If $r$ is odd, there is no isometry between $T(Y)_{\mathbb {Q}}$ and $T(X)_{\mathbb {Q}}$.
\item If $r$ is even, there exists an isometry between $T(Y)_{\mathbb {Q}}$ and $T(X)_{\mathbb {Q}}$ if and only if $T(Y)_{\mathbb {Q}}$ is isometric to $T(Y)_{\mathbb {Q}}(p)$. This property is equivalent to the following:

a) If $p=2$, for every prime number $q$ congruent to $3$ or $5$ modulo $8$, the $q$-adic valuation $\nu_q(\det T(Y))$ is even.

b) If $p>2$, for every prime number $q>2$, $q\neq p$, such that $p$ is not a square in $\mathbb F_q
$, the number $\nu_q(det(T_Y))$ is even and the following equation holds in $\mathbb F_p^*/(\mathbb F
_p^*)^2$
\[res_p(det(T_y)=(-1)^{\frac{p(p-1)}{2}+\nu_p(det(T_Y))}\]
where $res_p(det(T_Y))=\frac{det(T_Y)}{p^{\nu_p(det(T_Y))}}$.

\end{enumerate}
\end{theorem}
This theorem allows us to find $2$-isogenies as self isogenies on $Y_2$ and $Y_{10}$. 
In a previous paper we gave all the $2$-isogenies of $Y_2$ and exhibited self isogenies on $Y_2$.
In section $3$ we also exhibited $2$-isogenies as self isogenies on $Y_{10}$.
\begin{theorem}
All the $3$-isogenies between $Y_2$ with transcendental lattice $T(Y_2)=[2 \quad 0 \quad  4]$ and $Y_{10}$ with transcendental lattice  $T(Y_{10})=[6 \quad 0 \quad  12]$, all the $3$-isogenies between $Y_{10}$ and the $K3$ surfaces with transcendental lattices $[4 \quad 0 \quad 18]$ and $[2 \quad 0 \quad 36]$ are isometries of their rational transcendental lattices.

\end{theorem}  
In sections $5$ and $6$,
we proved that all the $3$- isogenies on $Y_2$ are between $Y_2$ and $Y_{10}$ and some $3$-isogenies on $Y_{10}$ are between $Y_{10}$ and other $K3$ surfaces with dicriminant $72$, namely $[4 \quad 0 \quad 18]$ or $[2 \quad 0 \quad 36 ]$. These results illustrate Boissi\`ere, Sarti and Veniani's theorem. Indeed $det(T(Y_2))=8$, hence $res_3(8)=8$ which is congruent modulo $3$ to $(-1)^3$ and $det(T(Y_{10}))=8\times 9$, hence $res_3(8\times 9)=8$ which is congruent modulo $3$ to $(-1)^{3+2}$. And, since
\[T_{\mathbb Q}(Y_2)=\begin{pmatrix}
	2 & 0\\
	0 & 1\\
\end{pmatrix}\qquad T_{\mathbb Q}(Y_{10})=\begin{pmatrix}
	6& 0\\
	0 & 3\\
      \end{pmatrix} \]
 \[   T_{\mathbb Q}([4 \quad 0 \quad 18])=\begin{pmatrix}
	1& 0\\
	0 & 2\\
\end{pmatrix}\qquad T_{\mathbb Q}([2 \quad 0 \quad 36])=\begin{pmatrix}
	2& 0\\
	0 & 1\\
\end{pmatrix} \]

we find, as expected, these matrices are isometric since 
\[\begin{pmatrix}
	6 & 0\\
	0 & 3\\
\end{pmatrix}=\begin{pmatrix}
	2 &1\\
	-1 & 1\\
\end{pmatrix}\begin{pmatrix}
	1 & 0\\
	0 & 2\\
\end{pmatrix}\begin{pmatrix}
	2 & -1\\
	1 & 1\\
      \end{pmatrix}\]
    
 \[   \begin{pmatrix}
	2& 0\\
	0 & 1\\
\end{pmatrix}=\begin{pmatrix}
	0 &1\\
	1 & 0\\
\end{pmatrix}\begin{pmatrix}
	1 & 0\\
	0 & 2\\
\end{pmatrix}\begin{pmatrix}
	0& 1\\
	1 & 0\\
\end{pmatrix}.\]

{\bf{Some remarks}}

As a consequence of Boissi\`ere, Sarti and Veniani's theorem, there could be $2$ or $3$-self isogenies on elliptic fibrations of $Y_2$ and $Y_{10}$. 
Indeed we found $2$-self isogenies on both $Y_2$ and $Y_{10}$. As for $3$-isogenies, there is no self-isogeny on $Y_2$ and also probably none on $Y_{10}$.

Concerning rank $0$ elliptic fibrations, using Shimada-Zhang's table \cite{SZ}, we recover easily all our results without using Weierstrass equations. We have only to know the transformation by a $2$- or a $3$-isogeny of a type of singular fiber. This can be obtained using Tate's algorithm \cite{T} and an analog of Dockchitser's remark \cite{Do}.





\begin{thebibliography}{25}





  

  
\bibitem{Be}
M. J. Bertin, {\it Mesure de Mahler et s\'erie $L$ d'une surface $K3$-singuli\`ere}, Publ. Math. Besan{\text{\c{c}}}on, Actes de la Conf\'erence Alg\`ebre et Th\'eorie des Nombres (2010), 5--28.



\bibitem{BL1}
M.-J. Bertin, and O. Lecacheux,  \textit{Elliptic fibrations on the modular surface associated to $\Gamma_1 (8)$}, in Arithmetic and geometry of K3 surfaces and Calabi-Yau threefolds,  Fields Inst.Commun., 67, Springer, New York, (2013), 153--199.


\bibitem{BL2} M.-J. Bertin, and O. Lecacheux,  \textit{Ap\'ery-Fermi pencil of $K3$-surfaces and 2-isogenies} J. Math. Soc. Japan (2) {\bf{72}}  (2020), 599--637.
  
\bibitem{BL3} M.-J. Bertin, and O. Lecacheux, \emph{Elliptic fibrations of a certain $K3$ surface of the Ap\'ery-Fermi pencil}, Publ. Math. Besan{\text{\c{c}}}on, Alg\`ebre et Th\'eorie des Nombres, 44, (2022), 5--36. 
  
\bibitem{BSV} S. Boissi\`ere, A. Sarti, D. Veniani, \emph{ On prime degree isogenies between $K3$ surfaces,} Rendiconti del Circolo Matematico di Palermo Series 2, Volume 66, Issue 1, (2017),  3--18.

\bibitem{Bour}
N. Bourbaki, \textit{Groupes et alg\` ebres de Lie}, Chap.4, 5, 6, (Masson, 1981).

  

\bibitem{Do}
T. and V. Dockchitser, \emph{A remark on Tate's algorithm and Kodaira types}, Acta Arith. 160 (2013), 95--100.


\bibitem{BE}
S. Brandhorst and N. Elkies, \emph{Equations for a K3 Lehmer map}, arXiv:2103.15101v2 [math.AG] 26 Sept 2022, 1--28.  to appear in J. Algebraic Geom. 


\bibitem{GS}
B. van Geemen and M. Sch\"{u}tt, \emph{On families of $K3$ surfaces with real multiplication}, arXiv:2310.05196v2 [math.AG] 24 Nov. 2023, 1--25.
  


\bibitem{GPM}
  A. Garbagnati, and Y. Prieto Monta\~ {n}ez, \emph{Generalized Shioda-Inose structure of order $3$}, Advances in Geometry 24, (2024) 183-207,
  arXiv:2209.10141v1 [math.AG] 21 Sep 2022, 1--32.  

\bibitem{Kum}
A. Kumar, \emph{Elliptic fibrations on a generic Jacobian Kummer surface}, J. Algebraic Geom. 23 (2014), 4, 599--667.  




\bibitem{KK}
A. Kumar, and M. Kuwata,  \textit{Elliptic }$K3$ \textit{ surfaces associated
with the product of two elliptic curves: Mordell-Weil lattices and their
fields of definition }, Nagoya Math. J. {\textbf{228}} (2017), 124--185. 


\bibitem{Ku1}
M. Kuwata,  \textit{Elliptic $K3$ surfaces with given Mordell-Weil rank}, Comment. Math. Univ. St. Paul.
{\textbf{49}} (2000), 91--100.



\bibitem{KU}  M. Kuwata, and  K. Utsumi, \emph{ Mordell-Weil lattice of Inose's elliptic K3 surface arising from the product of 3-isogenous elliptic curves,}
\ J. Number Theory \textbf{190} (2018), 333--351. 

\bibitem{L} O. Lecacheux, \emph{Weierstrass equations for all elliptic fibrations on the modular $K3$ surface associated to $\Gamma_1(7)$}, Rocky Mountain J. Math. 45 (5) (2015), 1481--1509.



  
  


\bibitem{Nik}
 V. V.  Nikulin, \textit{Integral symmetric bilinear forms and some of their applications}, Izv. Math. (1) {\textbf{14}},  (1980), 103--167.




 
 


\bibitem{Shim1}
I. Shimada,  \textit{$K3$ surfaces and Lattice Theory, Talk Dec. 2, 2014, Hano\"{i} Vietnam Institute for Advavanced Study in Mathematics, (Slides on Shimada's web page)}.

\bibitem{Shim2}
I. Shimada, \textit{On elliptic $K3$ surfaces}, arXiv:math 0505140v3 [math.AG] 11Apr 2006, 1--114.  

\bibitem{SZ}
I. Shimada,  and D. Q. Zhang, \textit{Classification of extremal elliptic $K3$ surfaces and fundamental groups of open $K3$ surfaces}, Nagoya Math. J. {\textbf{161}} (2001), 23--54.

\bibitem{Sh0}
  T. Shioda,   \textit{Kummer sandwich theorem of certain elliptic $K3$ surfaces }, Proc. Japan Acad., {\textbf{82}} Ser. A (2006), 137--140.
  
\bibitem{Sh1}
T. Shioda,  $K3$ \textit{surface and Sphere packing }, J.
Math. Soc. Japan {\textbf{60}} (2008), 1083--1105.



\bibitem{Sh2} T. Shioda,  \textit{Correspondence of elliptic curves and Mordell-Weil
lattices of certain elliptic K3 surfaces,}  Algebraic Cycles and Motives,
 Cambridge Univ. Press {\bf{2}} (2007) 319--339. 

\bibitem{Sh3} T. Shioda, \textit{On elliptic modular surfaces}, J. Math. Soc. Japan, 24 (1972), 20-59.
 
\bibitem{SI}
T. Shioda \& H. Inose, \emph{On singular $K3$ surfaces} in
 Complex analysis and algebraic geometry (Baily, Shioda T. ed.) , Cambridge (1977), 119-135.


 

  
 


\bibitem{T}
  J. Tate, \textit{Algorithm for determining the type of a singular fiber in an elliptic pencil}, in  Modular Functions of One Variable IV, Lect. Notes in Math. 476, B. J. Birch and W. Kuyk, eds, Springe
  r-Verlag, Berlin, 1975, 33-52.



\end{thebibliography}
\end{document}